\documentclass[11pt,a4paper]{article}
\usepackage{amsmath}
\usepackage{amssymb}
\usepackage{isolatin1}
\usepackage{txfonts}
\usepackage[french]{babel}

\newtheorem{definition}{Définition}[section]
\newtheorem{lem}[definition]{Lemme}
\newtheorem{theorem}{Théorème}
\numberwithin{equation}{section}

\newenvironment{dem}{\smallskip\noindent{\bf Démonstration}~:%
{\nopagebreak[0]}}%
{\nopagebreak[0]\hfill$\Box$ \smallskip}%

\newenvironment{rem}{\smallskip\noindent {\bf Remarque}~: }%
{\par\hspace*{7cm}\smallskip\par}

\newcommand{\al}{\alpha}

\newcommand{\lb}{\lambda}
\newcommand{\ro}{\rho}
\newcommand{\rok}{\rho_k}
\newcommand{\veps}{\varepsilon}
\newcommand{\intpart}[1]{\left\lfloor #1 \right\rfloor}

\renewcommand{\bar}[1]{\underline{#1}}
\newcommand{\bra}{\bar\alpha}

\newcommand{\brx}{\bar x}

\newcommand{\sg}{''}

\newcommand{\N}{\mbox{$\mathbb{N}$}}
\newcommand{\R}{\mbox{$\mathbb{R}$}}
\newcommand{\Rplus}{\R_{+}}

\newcommand{\ca}{{\cal A}}
\newcommand{\da}{\ca}
\newcommand{\cd}{{\cal D}}
\newcommand{\Primes}{{\cal P}}

\renewcommand{\AA}{{\rm\bf A}}
\newcommand{\BB}{{\rm\bf B}}

\newcommand{\kalm}{Kalm\'ar}

\newcommand{\dv}{\mid}
\newcommand{\st}{^\star}

\def\build#1_#2^#3{\mathrel{\mathop{\kern 0 pt#1}\limits_{#2}^{#3}}}

\newcommand{\abs}[1]{\left\lvert #1 \right\rvert}
\newcommand{\bigo}[1]{O\left(#1\right)}
\newcommand{\floor}[1]{\left\lfloor #1 \right\rfloor}
\newcommand{\nsup}[1]{\left\|#1\right\|}
\newcommand{\set}[1]{\left\{#1\right\}}
\newcommand{\setof}[2]{\left\{#1\ ;\ #2\right\}}

\newcommand{\dr}[1]{\frac{\partial}{\partial #1}}
\newcommand{\dron}[2]{\frac{\partial #1}{\partial #2}} 
\newcommand{\ddron}[2]{\frac{\partial^2 #1}{\partial #2^2}}
\newcommand{\dddron}[3]{\frac{\partial^2 #1}{\partial #2 \partial #3}}
\newcommand{\iem}{^{\text{ème}}}

\newcommand{\dsm}[1]{\mbox{$\displaystyle #1 $}}
\newcommand{\suppress}[1]{}

\title{Grandes valeurs et nombres champions de la fonction  
arithmétique de \kalm}
\author{M. Deléglise, M. O. Hernane, J.-L. Nicolas}
\date{}

\begin{document}
\maketitle

\renewcommand{\abstractname}{Abstract}
\newcommand{\omgb}{\varpiup}
\newcommand{\das}{\da^\star}

\begin{abstract}
The Kalm\'ar function $K(n)$ counts the factorizations
{$n=x_1 x_2 \dots x_r$} with $x_i \ge 2$, $(1 \le i \le r)$.  Its
Dirichlet series is $\sum_{n=1}^{\infty} \frac{K(n)}{n^s} =
\frac{1}{2-\zeta(s)}$ where $\zeta(s)$ denotes the Riemann $\zeta$
function. Let $\ro = 1.728\dots$ the root greater than $1$ of the
equation $\zeta(s) = 2$.  Improving on preceding results of Kalm\'ar,
Hille, Erd\"os, Evans, and Klazar and Luca, we show that there exist
two constants $C_5$ and $C_6$ such that, for all $n$, $K(n) \le \ro\log
n -C_5 (\log n)^{\frac{1}{\ro}}/\log\log n$ holds, while, for
infinitely many $n'$ s, we have $K(n) \ge \ro\log n -C_6 (\log
n)^{\frac{1}{\ro}}/\log\log n$.

An integer $N$ is called a $K$-champion number if $M < N \implies K(M)
< K(N)$. Several properties of $K$-champion numbers are given, mainly
about the size of the exponents and the number of prime factors in the
standard factorization into primes of a large enough $K$-champion
number.

The proof of these results is based on the asymptotic formula of
$K(n)$ given by Evans, and on the solution of a problem
of optimization.
\end{abstract}

\noindent{\bf Keywords~: }  Kalm\'ar's function, factorisatio
numerorum, highly composite numbers, champion numbers, optimization.
\medskip

\noindent{\bf 2000 Mathematics Subject Classification~:} 11A25, 11N37,
49K10.
\medskip

\section{Introduction}
Soit $\tau_r(n)$ le nombre de solutions de l'équation diophantienne
\begin{equation}\label{taur}
x_1 x_2 \dots x_r = n.
\end{equation}
On a
\[
\tau_0(n) = 
\begin{cases}
1 \text{ si } n=1 \\
0 \text{ si } n \ge 2,
\end{cases}
\quad \tau_1(n) = 1,
\quad \tau_2(n) = \sum_{d \dv n} 1,
\quad \tau_r(n) = \sum_{d \dv n} \tau_{r-1}(d).
\]
La série génératrice est
\[
\sum_{n=1}^{\infty} \frac{\tau_r(n)}{n^s} = \zeta(s)^r
\]
où $\zeta(s) = \sum_{n=1}^{\infty} \frac{1}{n^s}$ est la fonction
de Riemann.
\medskip

La fonction de \kalm\, ou \og factorisatio numerorum \fg\ compte le
nombre de solutions de \eqref{taur} pour tout $r$, mais
avec la restriction que chaque facteur $x_i$ doit vérifier $x_i \ge
2$. Ainsi les factorisations $12 = 6\times 2 = 4
\times 3 = 3 \times 4 = 3 \times 2 \times 2 = 2 \times 6 = 2 \times 3
\times 2 = 2 \times 2 \times 3$ donnent $K(12) = 8$.  On pose $K(1) =
1$. Pour $n \ge 2$, la fonction de \kalm\ satisfait
\begin{equation}\label{Kconv}
K(n) = \sum_{d \dv n,\,\, d \ge 2} K\left(\frac{n}{d}\right)
= \frac{1}{2}\sum_{d \dv n} K\left(\frac{n}{d}\right).
\end{equation}
La série de Dirichlet est
\begin{equation}\label{Kdir}
\sum_{n=1}^{\infty} \frac{K(n)}{n^s} = \frac{1}{2-\zeta(s)}\cdot
\end{equation}
Elle est reliée aux fonctions $\tau_r$ par la formule
\begin{equation}\label{Ktau}
K(n)  = \frac{1}{2} \sum_{r=0}^{\infty} \frac{\tau_r(n)}{2^r}\cdot
\end{equation}
La fonction $K(n)$ a été introduite par L. \kalm\  en
1931, dans \cite{Kalm1} et \cite{Kalm2}, où il montre que,
lorsque $x \to \infty$,
\[
\sum_{n \le x} K(n) \sim \frac{-1}{\ro\zeta'(\ro)}x^{\ro}
\]
($\ro = 1.728\dots$ est la racine positive de $\zeta(\ro) = 2$) et
donne une majoration du reste. Ceci a été précisé par Ikehara \cite{Ikeha}
et H.-K. Hwang \cite{Hwa}.

La majoration très simple $K(n) \le n^{\ro},\ (n \ge 1)$, obtenue
dans \cite{chor} et \cite{Coppers} a été récemmement améliorée par Klazar
et Luca \cite{luca} qui ont démontré $K(n) \le n^{\ro}/2,\ (n \ge 2)$,
en utilisant l'inégalité $K(nn') \ge 2K(n)K(n')$ vraie pour
tout $(n,n')$ satisfaisant $2 \le n \le n'$. On trouvera d'autres
informations sur la fonction de Kalm\'ar dans \cite{luca}, paragraphe 5.

Nous nous proposons dans cet article d'étudier les grandes valeurs de
la fonction $K(n)$. Ce sujet a déja été abordé par \kalm\
\cite{Kalm1}, \cite{Kalm2}, Erd\"os \cite{erdos}, Hille \cite{hill},
Evans \cite{Evans} (Th. 6 et 7) et Klazar et Luca \cite{luca}.

Soit $f$ une fonction arithmétique réelle ; appelons $f$-champion un
nombre $N$ tel que $n < N \implies f(n) < f(N)$. Les nombre
$\tau_2$-champions ont été appelés ``highly composite'' par Ramanujan
qui les a étudiés dans sa thèse \cite{rama}.

La fonction $O(n)$ d'Oppenheim (cf. \cite{macm1}, \cite{oppenheim},
\cite{canerpom}, \cite{HarrisSubbarao}) a la même définition que celle
de \kalm\  mais cette fois l'ordre des facteurs ne compte pas ; $12$ n'a
plus que les factorisations $12 = 6 \times 2 = 4 \times 3 = 3 \times
2\times 2$ et $O(12) = 4$.  Ainsi $O(n)$ compte le nombre de
\emph{partitions multiplicatives} de $n$ en parts $\ge 2$. Les nombres
$O$-champions, appelés ``highly factorable'' ont été étudiés dans
\cite{canerpom} et \cite{kim}.  A la fin de l'article \cite{canerpom},
le problème 5 demande quel est l'ordre maximum de la fonction $K(n)$,
et à quoi ressemblent les champions de $K(n)$. Les théorèmes
\ref{thomega}, \ref{decNth1}, \ref{decNth2}, apportent des éléments de
solutions à ce problème.

Soit $\ca \subset \set{2,3,4,\dots }$ ; dans \cite{hill} et
\cite{erdos} (cf.  aussi \cite{KnpKnpW1} et \cite{KnpKnpW2}) Hille et
Erd\"os ont généralisé la fonction de \kalm\  en définissant la fonction
$K_{\ca}(n)$ qui compte le nombre de solutions de \eqref{taur} pour
tout $r$ avec la restriction que chaque $x_i$ doit appartenir à $\ca$.

Dans l'article \cite{Her-Nic}, sont étudiées les grandes valeurs de la
fonction $K_{\Primes}(n)$, où $\Primes = \set{2,3,5,7,11,\dots}$ est
l'ensemble des nombres premiers, et quelques propriétés des nombres
$K_{\Primes}$-champions. Il est facile de voir que si la décom\-position
de $n$ en facteurs premiers est $n = q_1^{\al_1}\dots q_k^{\al_k}$
alors
\begin{equation}\label{KP}
K_{\Primes}(n) = 
\binom{\al_1 + \al_2 + \dots +\al_k}{\al_1,\al_2,\dots,\al_k}
= \frac{(\al_1+\al_2+\dots+\al_k)!}{\al_1! \al_2!\dots \al_k!}\cdot
\end{equation}

La formule \eqref{KP} ne s'étend pas à la fonction de \kalm. La
formule ci-dessous est due à Mac-Mahon \cite{macm2}, n$^{\circ}$ 80,
(cf. aussi \cite{KM}, formule (4))
\[
K(q_1^{\al_1} q_2^{\al_2}\dots q_k^{\al_k})
= \sum_{j=1}^{\al_1 + \al_2 + \dots + \al_k}
\sum_{i=0}^{j-1}(-1)^i\binom{j}{i} 
\prod_{h =1}^{k} \binom{\al_h+j-i-1}{\al_h},
\]
mais elle ne permet pas d'étudier les grandes valeurs de $K(n)$.
Cependant, à partir de \eqref{Ktau}, Evans a donné dans \cite{Evans}
une très jolie formule asymp\-totique pour $K(q_1^{\al_1}
q_2^{\al_2}\dots q_k^{\al_k})$ lorsque $\Omega(n) = \al_1 + \al_2
+\dots + \al_k$ tend vers l'infini. C'est à partir de cette formule
asymptotique que nous obtiendrons tous nos résultats.

Soit $\lb > 1$. L'article de Evans \cite{Evans} (cf. aussi \cite{Gro})
considère une fonction $K_{\lb}(n)$ plus générale dont la
série génératrice est
\[
\frac{\lb-1}{\lb-\zeta(s)} = 
\sum_{n=1}^{\infty}\frac{K_{\lb}(n)}{n^s} \cdot
\]

Lorsque $\lb = 2$, $K_2(n)$ est la fonction de \kalm\  $K(n)$. Les
nombres $K_{\lb}(n)$ sont les \emph{nombres eulériens généralisés}.
Le nombre eulérien $A(n,k)$ qui compte le nombre de permutations
de $n$ objets avec $k$ montées est relié aux nombres eulériens
généralisés par la formule
\[
\sum_{k=0}^{n} A(n,k)\lb^{k} = K_{\lb}(q_1 q_2 \dots q_n)
\]
où $q_1,q_2,\dots,q_n$ sont des nombres premiers distincts
(cf. \cite{Car}).

Les résultats de cet article pourraient s'étendre en remplaçant
$K(n)$ par $K_{\lb}(n)$ ; dans un souci de clarté nous nous sommes
limités à $\lb = 2$.

Dans le paragraphe 2, nous rappelons la formule asymptotique d'Evans,
et nous donnons quelques propriétés qui nous seront utiles par
la suite.

Dans le paragraphe 3, le théorème 2 donne un encadrement de $K(n)$
à l'aide de la fonction $F$. Le problème d'optimisation
\eqref{opt2} a été résolu par Evans \cite{Evans}, lemme 6, à l'aide
des multiplicateurs de Lagrange. Mais, afin de préciser
les comportement de $F$ au voisinage du maximum, nous déterminons
la forme quadratique des dérivées secondes, dont le calcul
présente, curieusement, des simplifications exceptionnelles. 

Au paragraphe 4, le théorème 3 améliore les théorèmes
3.1 et 4.1 de \cite{luca} et précise l'ordre maximum du
logarithme de la fonction de Kalm\'ar.

Les propriétés des nombres $K$-champions sont données
au paragraphe 5 par les théorèmes 4, 5, 6 et 7. Une table
de ces nombres figure en annexe. 

La démonstration des théorèmes 3, 4, 5, 6 et 7 suit d'assez
près la preuve des théorèmes correspondants de
\cite{Her-Nic}. Cependant, le remplacement de la formule
exacte (5) pour la fonction $K_{\cal P}$ par la formule
asymptotique de Evans (théorème 1) pour la fonction
$K$ complique les démonstrations. Ceci est particulièrement
net dans la preuve du théorème 6.

Nous avons plaisir à remercier L. Rifford pour l'aide
apportée à la résolution du problème d'optimisation
étudié au paragraphe 3.

\paragraph{Notations}
On utilisera les notations suivantes.
\begin{enumerate}
\item
Pour toute suite $\bar x = (x_i)_{1 \le i \le \omgb}$ de réels
de longueur $\omgb$, finie ou infinie, on note 
$\Omega(\brx) = \sum_{i=1}^{\omgb} x_i$
(lorsque cette somme a un sens)
et $\nsup{\brx} = \sum_{i=1}^{\omgb} \abs{x_i}$.
Si $\brx=(x_1,x_2,\dots,x_k) \in \Rplus^k$ et
${\underline y} = (y_1,y_2,\dots,y_\ell) \in \Rplus^\ell$,
on définit $\brx'=(x_1',x_2',\dots)$ par
$x_i'=x_i$ pour $1 \le i \le k$, et $x_i'=0$ pour $i > k$
et de même ${\underline y'} = (y_1',y_2',\dots)$
par $y_i' = y_i$ pour $i \le \ell$ et $y_i'=0$
pour $i > \ell$. Par définition on pose
$\nsup{\brx -{\underline y}} = \nsup{\brx'-\underline{y}'}$.

\item On note $\da$ l'ensemble des suites de réels positifs ou nuls
telle que $0 \le \Omega(\brx) < +\infty$, 
$\underline 0 = (0,0,0,\dots)$, 
et $\das = \da \setminus \set{\underline 0}$
\item Si $\bar x \in \das$ on note 
$\omgb(\bar x) = \sup\setof{j\in \N}{x_j \ne 0}$.
\item
$p_k$ représente le $k\iem$ nombre premier. Par le théorème
des nombres premiers on sait que $p_k\sim k\log k$
lorsque $k \to \infty$.
\item
Pour $n$ entier, de décomposition en facteurs premiers
$n = q_1^{\al_1}q_2^{\al_2}\dots q_k^{\al_k}$, on note
\[
\omega(n) = k 
\quad\text{ et} \quad
\Omega(n) = \sum_{i=1}^{k} \al_i.
\]
\suppress{
\item
Pour tout $x = (x_0,x_1,x_2,\dots) \in \R^{\N}$ on note
\[
\nsup{x} = \sum_{j=0}^{\infty} \abs{x_j}
\]
}
\end{enumerate}

\newcommand{\cxi}{c(\bar{x})+x_i}
\newcommand{\cxj}{c(\bar{x})+x_j}
\newcommand{\cxii}{c(\bar{\xi})+\xi_i}
\newcommand{\str}{^{\star}}

\section{L'estimation de Evans}

\subsection{La fonction $c$}

La fonction $c(\brx)$ définie ci-dessous a été introduite par Evans
dans \cite{Evans} lorsque $\brx$ est de longueur
$\omgb$ finie. Nous l'étendrons ici aux suites infinies.

\begin{definition}\label{defc}
Soit $\brx \in \das$ et $\omgb = \omgb(\brx)$.
Il existe un unique $c = c(\bar x)> 0$ tel que
\begin{equation}\label{defceq}
\prod_{j=1}^{+\omgb} \left(1+\frac{x_j}{c} \right) 
= 2.
\end{equation}
De plus, pour tout $\lb > 0$,
\begin{equation}\label{chomogen}
c(\lb x_1, \lb x_2,\dots) = \lb c(x_1,x_2,\dots)
\end{equation}
et
\begin{equation}
\Omega(\bar x) \le c(\bar x) 
\le \frac{\Omega(\bar x)}{\log 2} \le \frac{3\Omega(\bar x)}{2}\cdot
\label{defcenc}
\end{equation}
\end{definition}

\begin{dem}
Pour tout $t > 0$ on note
\[
H(\bar x,t) = \sum_{j=1}^{\omgb} \log\left(1+\frac{x_j}{t}\right)
\le \frac{1}{t}\nsup{\brx}.
\]
La fonction $t \mapsto H(\bar x,t)$ décroît de $+\infty$ à
$0$ lorsque $t$ croît de $0$ à $+\infty$.  Ceci assure
l'existence et l'unicité de $c$.
La propriété \eqref{chomogen} est immédiate. Prouvons l'encadrement
\eqref{defcenc}.
Pour $t=\Omega(\bar x) = \Omega$,
\[
H(\bar x,\Omega) = 
\log\prod_{j=1}^{\omgb} \left(1+\frac{x_j}{\Omega}\right)
\ge \log\left(1+\sum_{j=1}^{\omgb} \frac{x_j}{\Omega}\right) =\log  2.
\]
Ceci montre que $c \ge \Omega$.
La majoration de $c(\bar x)$ résulte de 
\[
\log 2 = \sum_{j=1}^{\omgb} \log\left(1+\frac{x_j}{c}\right)
\le \sum_{j=1}^{\omgb} \frac{x_j}{c} = \frac{\Omega}{c}\cdot
\]
\end{dem}

\begin{rem}
  On vérifie immédiatement que si $\omgb$ est fini on a
\begin{equation}\label{cxk0}
c(x_1,x_2,\dots,x_\omgb,0)=  c(x_1,x_2,\dots,x_\omgb)  
\end{equation}
et que la suite infinie $\bar x'$ définie par $x_i'=x_i$ pour $i
  \le \omgb$ et $x_i' = 0$ pour $i > \omgb$ est un élément de
  $\da$ et $c(\brx') = c(x_1,x_2,\dots,x_{\omgb})$. Notons aussi que
  $c$ est symétrique en les $x_i$, et que c'est une fonction
  croissante de chaque variable $x_i$. Par convention, on pose
  $c(\underline 0) = 0$.
\end{rem}

\begin{lem}\label{cktoc}
Soit $\bar x = (x_1,x_2,\dots) \in \das$, $c = c(\bar x)$ et,
pour tout $k \ge 1$
\[
c_k = c(x_1,x_2,\dots,x_k,0,0,\dots).
\]
Alors 
\[
\lim_{k \to \infty} c_k = c.
\]
\end{lem}

\begin{dem}
La suite $c_k$, croissante et majorée par $c$ $\Big($car
$\prod_{j=1}^{k}\left(1+\frac{x_j}{c}\right) < 2$$\Big)$ admet une limite
$\hat c$. Comme $\prod_{i=1}^{k} \left(1+\frac{x_i}{t}\right)$
converge uniformément vers $\prod_{i=1}^{+\infty}
\left(1+\frac{x_i}{t}\right)$, sur tout intervalle $[u,+\infty[$,
$u>0$, on a
\[
\prod_{i=1}^{+\infty} \left(1+\frac{x_i}{\hat c}\right)
= \lim_{k \to +\infty}
\prod_{i=1}^{k} \left(1+\frac{x_i}{c_k}\right)
= 2 =
\prod_{i=1}^{+\infty} \left(1+\frac{x_i}{c}\right),
\]
et cela donne $\hat c = c$.
\end{dem}

\begin{lem}\label{derivc}
Soit $\brx \in \das$.
Pour tout entier $i \ge 1$, $c$ admet une dérivée
partielle par rapport à $x_i$, qui est donnée par
\begin{equation}\label{dronc}
\dron{c(\bar{x})}{x_i} = \frac{1}{T(\bar{x})} \frac{c(\bar x)}{\cxi},
\quad\text{ avec }\quad
T(\bar{x}) = \sum_{i=1}^{\infty}  \frac{x_i}{\cxi}\cdot
\end{equation}
\end{lem}

\begin{dem}
\suppress{
  Puisque $c$ est une fonction symétrique de ses arguments on peut
  supposer $i=1$.  Pour $x_2,x_3,\dots$ fixés, considérons la
  fonction
\[
H(u,t) = \log\left(1+\frac{u}{t}\right)
+ \sum_{i=2}^{+\infty} \log\left(1+\frac{x_i}{t}\right)
- \log 2.
\]
Par défintion de $c$, $H(x_1,c(\bar x)) = 0$.
Soit $U$ l'ouvert de $\R^2$ formé des couples $(u,t)$ 
satisfaisant les inégalités
\[
\begin{cases}
t > \Omega/2\\
u > 0
\end{cases}
\text{ si } x_1 > 0
\quad\text{ et }\quad
\begin{cases}
t > \Omega/2\\
\abs{u} < t
\end{cases}
\text{ si } x_1 = 0.
\]
$U$ contient le point $(x_1,c(\bar x))$, et 
$H$ admet dans $U$ des dérivées partielles en $u$ et en $t$,
qui sont
\[
\dron{H}{u} = \frac{1}{u+t},
\quad
\dron{H}{t} = - \frac{\frac{u}{t^2}}{1+\frac{u}{t}}
-\sum_{i=2}^{+\infty} \frac{\frac{x_i}{t^2}}{1+\frac{x_i}{t}}
= -\frac{1}{t}\left( \frac{u}{u+t}
+\sum_{i=2}^{+\infty}\frac{x_i}{x_i+t}\right).
\]
L'existence de la dérivée partielle en $t$, provient de ce que
la série $\sum_{i=2}^{+\infty}\frac{x_i}{x_i+t}$ converge
uniformément en $t$, dans le domaine $t > \Omega/2$.
La dérivee partielle en $t$ est strictement positive en le point
$(x_1,c(\bar x))$. Par le théorème des fonctions implicites
la fonction $c$, définie pour $x_2,x_3,\dots$ fixés,
par l'équation $H(x_1,c) = \log 2$ est donc dérivable en
$x_1$, de dérivée
\[
\dron{c(\bar x)}{x_1} = 
-\,\,\dfrac{\dron{H(u,t)}{u}}
{\dron{H(u,t)}{t}}(x_1,c)
= 
\frac{\frac{c}{c+x_1}}{\sum_{i=1}^{+\infty}\frac{x_i}{x_i+c}}\cdot
\]
} Puisque $c$ est une fonction symétrique de ses arguments on peut
supposer $i=1$.  Fixons $x_2,x_3,\dots,x_k,\dots$, et supposons les
d'abord non tous nuls. Posons $\brx = (x_1,x_2,\dots,x_k,\dots).$
L'application $x_1 \mapsto c(\brx)$ est une bijection croissante de
$[0,+\infty[$ sur $[c_0,+\infty[$, avec $c_0 = c(0,x_2,x_3,\dots)> 0$,
car, par définition, $x_1$ s'explicite en fonction de $c$ par
\[
x_1 = \frac{2c}{\Pi}-c
\quad\text{avec}\quad
\Pi = \prod_{j=2}^{\infty}\left(1+\frac{x_j}{c} \right)\cdot
\]
La convergence de la série $\sum_{j=2}^{\infty} x_j$ entraine que
$\log\Pi = \sum_{j=2}^{\infty} \log\left(1+\frac{x_j}{c}\right)$ est une
fonction dérivable de $c$ sur $[c_0,+\infty[$ et que l'on a
\[
\frac{d\Pi}{\Pi dc} = -\frac{1}{c}\sum_{j=2}^{\infty}\frac{x_j}{c+x_j}
= -\frac{1}{c}\left(T(\brx)-\frac{x_1}{c+x_1} \right)\cdot
\]
Il en résulte que $c
\mapsto x_1$ est une bijection croissante et continûment dérivable de
$[c_0,+\infty[$ sur $[0,+\infty[$ et que l'on a,
puisque $\Pi = \dfrac{2c}{c+x_1}$,
\begin{equation}\label{dx1dc}
\frac{dx_1}{dc}= \frac{2}{\Pi}-\frac{2c}{\Pi^2}\frac{d\Pi}{dc}-1
= \frac{2}{\Pi}\left(1+T(\brx)-\frac{x_1}{c+x_1}\right)-1
= \left(\frac{c+x_1}{c} \right) T(\brx)\cdot
\end{equation}
Par le theorème d'inversion $c$ est une fonction continûment dérivable
de $x_1$ et \eqref{dx1dc} implique \eqref{dronc}.

Si $0=x_2=x_3=\dots$, la définition \eqref{defc} donne $c(\brx)=x_1$ et
le résultat est encore vrai.
\end{dem}

\begin{lem}\label{lemprod}
Soit $(\gamma_i)_{i \ge 1}$ une suite de réels
vérifiant $0 \le \gamma_i < 1$
et $\sum_{i=1}^{+\infty} \gamma_i < +\infty$.
Alors
\[
\prod_{i=1}^{+\infty} (1-\gamma_i) \ge 1 - \sum_{i=1}^{+\infty} \gamma_i.
\]
\end{lem}

\begin{dem}
Notons $S_n = \sum_{i=1}^{+\infty} \gamma_i^n$.
Le résultat est évident si $S_1  > 1$.
Supposons donc $S_1 \le 1$.
Alors, pour tout $n$, on a $S_n \le S_1^n$ et
\[
\log \prod _{i=1}^{+\infty}(1-\gamma_i)
= \sum_{i=1}^{+\infty} \log(1- \gamma_i)
= -\sum_{i=1}^{+\infty} \sum_{n=1}^{+\infty} \frac{\gamma_i^n}{n}
=  -\sum_{n=1}^{+\infty} \frac{S_n}{n}
\ge -\sum_{n=1}^{+\infty} \frac{S_1^n}{n}
\]
tandis que
\[
\log \bigg(1 - \sum_{i=1}^{+\infty} \gamma_i\bigg) 
= \log(1-S_1)
= -\sum_{n=1}^{+\infty} \frac{S_1^n}{n}\cdot
\]
\end{dem}

\begin{lem}\label{encadT}
  Pour tout $\bar x \in \das$, la quantité $T(\bar{x}) =
  \sum_{i=1}^{+\infty} \frac{x_i}{\cxi}$ vérifie
\[
\frac{1}{2} \le T(\bar x) \le 1.
\]
Et chaque dérivée partielle de $c$ vérifie
\[
0 \le \dron{c(\bar x)}{x_i} \le \frac{1}{T(\bar x)} \le 2.
\]
\end{lem}

\begin{dem}
L'encadrement \eqref{defcenc} donne la majoration
\[
T(\bar x) = \sum_{i=1}^{+\infty} \frac{x_i}{c+x_i} 
\le \sum_{i=1}^{+\infty} \frac{x_i}{c}
= \frac{\Omega}{c} \le 1.
\]
Pour la minoration, notons $\gamma_i = \frac{x_i}{c+x_i}$. Alors, par
la définition \ref{defc} de $c$,
\[
\prod_{i=1}^{\infty} (1-\gamma_i) = 
\prod_{i=1}^{\infty} \frac{c}{c+x_i} = \frac{1}{2}\cdot
\]
Le lemme \ref{lemprod} donne alors
\[
T(\brx) = \sum_{i=1}^{\infty} \gamma_i\ge 1 - \prod_{i=1}^{\infty} (1-\gamma_i)
= 1-\frac{1}{2} = \frac{1}{2}\cdot
\]
Pour le deuxième point, le lemme \ref{derivc} entraîne
\[
\dron{c(\bar{x})}{x_i} = \frac{1}{T(\bar{x})} \frac{c(\bar x)}{\cxi}
\le \frac{1}{T(\bar x)} \le 2.
\]
\end{dem}

\suppress{
\begin{lem}
  Soit $\bar x = (x_1,x_2,\dots) \in \da$ 
  avec $x_1 \ge x_2 \ge 1$. Soit $t$, $0 < t_0 \le x_2$,
  et $\bar y = (x_1 + 1,x_2-1,x_3,\dots)$.
\[
c(\bar y) < c(\bar x). 
\] 
\end{lem}

\begin{dem}
Soit $c = c(\bar x)$.
Posons $f(t) = c(x_1+t,x_2-t,x_3,\dots)$.
Alors, il existe $t$, $0 < t_0 < 1$, tel que
\[
c(\bar y) - c(\bar x) =  f'(t_0)
= \dron{c}{x_1}(\bar z) -  \dron{c}{x_2}(\bar z)
\]
avec $\bar z = (x_1+t_0,x_2-t_0,x_3,\dots,)$. 
Ensuite, le lemme \ref{derivc} donne
\[
 \dron{c}{x_1}(\bar z) -  \dron{c}{x_2}(\bar z) =
\frac{c(\bar z)}{T(\bar z)}
\left(\frac{1}{c(\bar z)+x_1+t_0)} - \frac{1}{c(\bar
    z)+x_2-t_0)}\right)
< 0.
\]
\end{dem}
}

\begin{lem}\label{lemcxy}
Soit $\bar x$ et $\bar x' \in \da$.
Alors
\[
\abs{c(\bar x') - c(\bar x)} \le 2\,\sum_{i=1}^{+\infty}\abs{x_i'-x_i}
= 2 \nsup{\brx'-\brx}.
\]
\end{lem}

\begin{dem}
\begin{enumerate}
\item
Si $\bar x'$ et $\bar x$ sont tous deux de même longueur finie
$k$ on note
\[
F(t) = c\big(\bar x + t(\bar x' - \bar x)\big).
\]
Alors $c(\bar x')-c(\bar x) = F(1)-F(0)$. De plus,
$F$ est dérivable sur l'intervalle $(0,1)$, de dérivée
\[
F'(t) 
= \sum_{i=1}^{k} (x_i'-x_i)\dron{c}{x_i}\big(\bar x + t(\bar x'-\bar x)\big).
\]
Par le théorème des accroissements finis et 
le lemme \ref{encadT}, il vient
\[
\abs{c(\brx)-c(\brx')} = \abs{F(1)-F(0)}
\le \sup_{0 < t < 1} \abs{F'(t)}
\le 2 \sum_{i=1}^{k}\abs{x_i'-x_i}.
\]
\item
Revenons maintenant au cas général.
Pour tout $k \ge 1$, notons
\[
c_k = c(x_1,x_2,\dots,x_k,0,0,\dots)
\text{ et }
c_k' = c(x_1',x_2',\dots,x_k',0,0,\dots)
\]
Par le premier point, on a
$
\abs{c_k' -c_k} \le 2\,\sum_{i=1}^{k} \abs{x_i'-x_i}
\le 2\,\nsup{\brx'-\brx}.
$
Avec le lemme \ref{cktoc} on en déduit
\[
\abs{c'-c} = \lim_{k \to \infty} \abs{c_k'-c_k}
\le 2\,\nsup{\brx'-\brx}.
\]
\end{enumerate}
\end{dem}

\begin{lem}\label{Tlipsch}
Soit $\brx,\brx' \in \da$.
Notons $\Omega = \Omega(\brx)$ et $\Omega' = \Omega(\brx')$.
Pour la fonction $T$ définie en \eqref{dronc}, on a la majoration
\[
\abs{T(\brx') - T(\brx)} \le 
\frac{3}{\max(\Omega,\Omega')} \nsup{\brx'-\brx}.
\] 
\end{lem}

\begin{dem}
  Notons $c=c(\bar x)$ et $c'=c(\bar x')$ et supposons
  $\Omega'\ge \Omega$.  Alors
\begin{eqnarray*}
T(\brx') - T(\brx) &=& \sum_{i=1}^{+\infty}
\left[\frac{x_i'}{c'+x_i'} - \frac{x_i}{c+x_i}\right]
=  \sum_{i=1}^{+\infty}
\frac{c x_i' - c' x_i}{(c'+x_i')(c+x_i)}\\
&=&
\sum_{i=1}^{\infty} \frac{c(x_i'-x_i)}{(c'+x_i')(c+x_i)} + 
\sum_{i=1}^{\infty} \frac{(c-c')x_i}{(c'+x_i')(c+x_i)}.
\end{eqnarray*}
et, en utilisant le lemme \ref{lemcxy} et \eqref{defcenc},
\begin{eqnarray*}
\abs{T(\brx') - T(\brx)} &\le& 
  \sum_{i=1}^{+\infty} \frac{\abs{x_i'-x_i}}{c'+x_i'}
  + \sum_{i=1}^{+\infty} \frac{x_i \abs{c'-c}}{(c'+x_i')(c+x_i)}\\ 
&\le&
  \sum_{i=1}^{+\infty} \frac{\abs{x_i'-x_i}}{\Omega'}
  + \abs{c'-c}\sum_{i=1}^{+\infty} \frac{x_i}{\Omega\Omega'}\\
&=&
  \frac{\nsup{\brx'-\brx}}{\Omega'}
  + \frac{\abs{c'-c}}{\Omega'} 
\le \frac{3}{\Omega'}\nsup{\brx'-\brx}.
\end{eqnarray*}
\end{dem}

\suppress{
Le lemme suivant affine la majoration du lemme \ref{lemcxy}.
\begin{lem}\label{lemcxy2}
Soit $x',x \in \da$, et 
\[
M = \sup_{t \in [0,1]} T\left(\bar x + t(\bar x'-\bar x)\right).
\]
Alors
\[
\abs{c(x')-c'(x)} \le \frac{1}{M} \nsup{x'-x}.
\]
\end{lem}

\begin{dem}
Notons
$X_k = (x_1,x_2,\dots,x_k,0,\dots)$,
$X_k' = (x_1',x_2',\dots,x_k',0,\dots)$ et
\[
M_k = \sup_{t \in [0,1]} T\left(X_k + t(X_k'-X_k)\right).
\]
Soit $F$ la fonction de $t$ définie sur $[0,1]$ par
\[
F(t) = c\big(X_k + t(X_k' - X_k)\big).
\]
Alors $c(X_k')-c(X_k) = F(1)-F(0)$. De plus,
$F$ est dérivable, de dérivée
\[
F'(t) 
= \sum_{i=1}^{k} (x_i'-x_i)\dron{c}{x_i}\big(X_k + t(X_k'-X_k)\big).
\]
Par définition de $M_k$, pour tout $t \in [0,1]$,
\[
\abs{F'(t)} \le \frac{1}{T\big(X_k + t(X_k'-X_k)\big)}
\sum_{i=1}^{k}\abs{x_i'-x_i}
\le \frac{1}{M_k} \sum_{i=1}^{k}\abs{x_i'-x_i}.
\]
Par le théorème des accroissements finis
\[
\abs{c(X_k')-c(X_k)} \le \frac{1}{M_k} \sum_{i=1}^{k}\abs{x_i'-x_i}
\le  \frac{1}{M_k} \nsup{\brx'-\brx}
\]
Par le lemme précédent
\end{dem}
}

\subsection{Approximation de $K(n)$}
Les définitions suivantes ont été introduites par Evans (cf. \cite{Evans}).

\begin{definition}\label{definitionalpha}
Pour $\brx = (x_1,x_2,\dots,x_{k}) \in \Rplus^{\star k}$ on définit
\begin{equation}\label{defalpha}
\AA(\brx) = \frac{1}{2\sqrt 2} e^{-\Omega(\brx)} 
\prod_{i=1}^{k} \frac{(\cxi)^{x_i}}{\Gamma(x_i+1)}\cdot
\end{equation}
\end{definition}

\begin{definition}
Pour $\brx = (x_1,x_2,\dots,x_{k}) \in \Rplus^{\star k}$ on définit
\begin{equation}\label{defB}
\BB(\brx) =  \left\{
\frac{1}{2c(\brx)}\sum_{i=1}^{k}\frac{x_i}{c(\brx) + x_i}
\right\}^{-\frac12}
= \sqrt{\frac{2c(\brx)}{T(\brx)}}\cdot
\end{equation}
\end{definition}
Du lemme \ref{encadT}, il résulte immédiatement le suivant

\begin{lem}
Pour $\brx  \in \Rplus^{\star k}$, on a
\begin{equation}\label{encb}
\sqrt{2c(\brx)} \le \BB(\brx) \le 2 \sqrt{c(\brx)}.
\end{equation}
\end{lem}
On trouvera également dans \cite{Evans} la démonstration du théorème
suivant.
\begin{theorem}\label{evans}[Evans]
Pour tout $\eta$, $0 \le \eta < 1/2$, il existe
$\Omega_0$ et $C_0$ tels que, pour tout
entier $n$ dont la décomposition en facteurs
premiers $n = q_1^{\al_1}q_2^{\al_2}\dots q_k^{\al_k}$
satisfait $\Omega(n) = \al_1 + \al_2 + \dots + \al_k \ge \Omega_0$,
on a, en posant $\bra = (\al_1,\al_2,\dots,\al_k)$,
\[
K(n)  = \sqrt\pi \AA(\bra) \BB(\bra) (1+ R(\bra))
\quad\text{ avec }\quad
\abs{R(\bra)} \le C_0\, (\Omega(\bra))^{-\eta}.
\]
\end{theorem}

\begin{rem}
La démonstration d'Evans est effective, et permettrait d'expli\-citer
des valeurs de $\Omega_0$ et $C_0$ pour une valeur donnée de $\eta$.
Le calcul est cependant technique et nous ne le ferons pas.
\end{rem}

\begin{lem}
Il existe deux constantes absolues $C_1$ et $C_2$ telles que,
pour tout entier $n\ge 2$, de décomposition en facteurs premiers
$n = q_1^{\al_1}q_2^{\al_2}\dots q_k^{\al_k}$, on ait
avec $\bra=(\al_1,\al_2,\dots,\al_k)$
\begin{equation}\label{ckc}
C_1 \sqrt\pi\, \AA(\bar \al) \BB(\bar \al) 
\le K(n) \le 
 C_2\sqrt\pi\, \AA(\bar \al)\BB(\bar \al).
\end{equation}
\end{lem}

\begin{dem}
Choisissons, par exemple $\eta = 1/4$, 
et $\veps > 0$ arbitraire; par le théorème \ref{evans} il existe
un $\Omega_1$ tel que
le rapport $K(n)/(\sqrt\pi\, \AA(\bar \al)\BB(\bar \al))$ 
diffère de $1$ de au plus $\veps$, pourvu que
$\Omega(n) = \Omega(\bar \al) \ge \Omega_1$.
Pour tous les entiers $n$ avec $1 \le \Omega(n) < \Omega_1$,
le nombre des valeurs de $\bar \al$ est fini.
\end{dem}

\begin{rem}
Nous avons calculé, 
pour chaque $r \in \set{1,2,\dots,20}$, et 
pour tous les $n$ tels que $\Omega(n) = r$,
le rapport $K(n)/(\sqrt\pi \AA(\bar \al)\BB(\bar \al))$.
Nous avons vérifié que ce rapport est minimum
lorsque $\bar \al = (r)$, c'est à dire lorsque $n$
est une puissance de nombre premier, $n=p^r$.
Il atteint son maximum lorsque $\bar \al = (1,1,\dots,1)$, c'est
à dire lorsque $n$ est un produit de $r$ facteurs premiers distincts.
Il est raisonnable de penser que cette propriété
est encore vraie pour toutes les valeurs de $r \ge 20$.
Ceci donne la conjecture suivante~: 
\emph{la valeur optimale de $C_2$ est
$1.084\, 437\, 552\dots$, atteinte pour $\bar \al = (1)$,
c'est à dire lorsque $n$ est premier. Et la valeur
optimale de $C_1$ est $1$, approchée
par les entiers sans facteurs carrés, lorsque le
nombre de leurs facteurs premiers tend vers l'infini.}
\end{rem}

\subsection{Les constantes  $\ro$, $\ro_k$, $a$ et $a_k$}
Les constantes $\ro,\rok,a,a_k$ ont été introduites par Hille
\cite{hill}, Evans  (\cite{Evans}, p. 169) et Klazar
et Luca \cite{luca}.
Dans ce paragraphe nous rappelons et précisons leur comportement.

\begin{definition}
Soit $\zeta(s) = \prod_{p}\left(1-\frac{1}{p^s}\right)^{-1}
= \sum_{n=1}^{+\infty} \frac{1}{n^s}$
la fonction de Riemann, et pour tout $k \ge 1$,
posons $\zeta_k(s) = \prod_{j=1}^{k}\left(1-\frac{1}{p_j^s}\right)^{-1}$.

\noindent
On définit $\ro = 1.728 647 238 998\dots$, et pour tout $k \ge 1$,
$\ro_k > 0$ par
\begin{equation}\label{defro}
\zeta(\ro) = 2,
\quad
\zeta_k(\ro_k) = 2.
\end{equation}
Pour $s > 1$, on définit $L(s) = -\log\zeta(s)$, de sorte que
$$
L'(s) = -\frac{d}{ds}(\log\zeta(s)) = 
\sum_{i=1}^{+\infty} \frac{\log p_i}{p_i^{s}-1},
$$
et  $L_k(s) = -\log\zeta_k(s)$.
On définit $a = 1.100020011\dots$ et $a_k$ par
\begin{equation}\label{defaak}
\frac{1}{a} = \sum_{i=1}^{+\infty} \frac{\log p_i}{p_i^{\ro}-1} 
= L'(\ro)
\quad\text{ et }\quad
\frac{1}{a_k} = \sum_{i=1}^{k} \frac{\log p_i}{p_i^{\ro_k}-1} 
= L_k'(\ro_k).
\end{equation}
\end{definition}
\medskip

\noindent
La figure \ref{figro} donne la valeur des premiers termes des 
suites $(\ro_k)$ et $(a_k)$.

\begin{lem}\label{romrokprop}
La suite $(\ro_k)_{k \ge 1}$ est croissante
et $\lim_{k \to +\infty} \ro_k = \ro$.
De plus,
lorsque $k$ tend vers l'infini,
\begin{equation}\label{romrok}
\rho -\rho_k \sim  \frac{2}{(-\zeta'(\rho))(\rho-1)k^{\rho-1}(\log k)^{\rho}}
= \frac{1.509\dots}{k^{\ro-1}(\log k)^{\ro}}\cdot
\end{equation}
\end{lem}

Ce lemme est démontré en \cite{luca},
paragraphe 2. En introduisant l'exponentielle
intégrale
$E_1(x) = \int_{x}^{\infty} \frac{e^{-t}}{t}dt$,
il est possible de démontrer, comme dans \cite{Her-Nic},
 (3.16)
\[
\ro-\ro_k = \frac{2 E_1((\ro-1)\log p_k)}{-\zeta'(\ro)}
+ \frac{{\cal O}_u(1)}{k^{\ro-1}(\log k)^{u}}
\]
pour tout $u > 0$.

\begin{lem}\label{akma}
Soient $a_k$ et $a$ définis en \eqref{defaak}.
La suite $(a_k)$ est décroissante et,
lorsque $k$ tend vers l'infini, on a l'équivalence
\[
a_k -a \sim \frac{a^2}{\ro-1}\frac{1}{(k\log k)^{\ro-1}} \cdot
\]
\end{lem}

\begin{figure}
\[
\begin{array}{|c|c|c|c|c|c|c|c|}
\hline
k    = & 1 & 2 & 3 & 10 & 100 & 1000  & \infty\\
\hline
\ro_k= &  1.00000& 1.43527& 1.56603& 1.69972 & 1.72658 
   & 1.72843      & 1.72864   \\
a_k = & 1.44269   &  1.44336   &  1.36287   & 1.19244  & 1.11279 
& 1.10196 & 1.1000200\\  
\hline
\end{array}
\]
\caption{Quelques valeurs de $\ro_k$ et $a_k$}
\label{figro}
\end{figure}

\begin{dem}
La décroissance de la suite $(a_k)$ résulte de la croissance
de $(\rok)$ et de \eqref{defaak}.  
Lorsque $k \to \infty$, il résulte de \eqref{defaak}
et de $\lim \ro_k = \ro$ (lemme \ref{romrokprop})
que $\lim a_k = a$ et  
\begin{equation}\label{invega}
\frac{1}{a}-\frac{1}{a_k} = \frac{a_k-a}{aa_k} 
\sim \frac{a_k-a}{a^2}\cdot
\end{equation}
Calculons un équivalent de $1/a-1/a_k$.
Les définitions \eqref{defaak} donnent
\begin{equation}\label{invega2}
\frac{1}{a} - \frac{1}{a_k} = L'(\ro)  - L_k'(\rok) 
= L'(\ro) - L_k'(\ro) + L_k'(\ro) - L_k'(\ro_k).
\end{equation}

\paragraph{Estimation de $L'(\ro) - L_k'(\ro)$}
On part de
\begin{equation}\label{lkml}
L'(\ro) - L_k'(\ro)  = \sum_{i=k+1}^{+\infty} 
\frac{\log  p_i}{p_i^{\ro}-1}\cdot
\end{equation}
Les équivalences (cf., par exemple, \cite{Dieu}, (3.11.10),
(3.11.6) et (3.10.5))
\begin{eqnarray*}
\sum_{i=k+1}^{+\infty}\frac{\log p_i}{p_i^{\ro}-1} &\sim&
\sum_{i=k+1}^{\infty}\frac{\log(i\log i)}{(i\log i)^{\ro}} \sim
\sum_{i=k+1}^{\infty}\frac{1}{i^{\ro}(\log i)^{\ro-1}} \sim
\int_{k}^{\infty}\frac{dt}{t^{\ro}(\log t)^{\ro-1}}\\
 &\sim&
\frac{1}{(\ro-1)}\,\frac{1}{k^{(\ro-1)}(\log k)^{\ro-1}}
\end{eqnarray*}
donnent avec \eqref{lkml}
\begin{equation}\label{lkmleq}
L_k'(\ro) - L'(\ro)  \sim 
\frac{1}{(\ro-1)}\,\frac{1}{k^{(\ro-1)}(\log k)^{\ro-1}}\cdot
\end{equation}

\paragraph{Estimation de $L_k'(\ro) - L_k'(\ro_k)$}
La définition de la dérivée et \eqref{romrok}
donnent
\begin{multline}\label{ts}
L_k'(\ro) - L_k'(\ro_k) \sim (\ro - \ro_k) L_k\sg(\ro)
\sim
-\frac{2 L_k\sg(\ro)}{(\ro-1)\zeta'(\ro)}\frac{1}{k^{\ro-1}(\log
 k)^{\ro}}\\
\sim
-\frac{2 L\sg(\ro)}{(\ro-1)\zeta'(\ro)}\frac{1}{k^{\ro-1}(\log
  k)^{\ro}}
\cdot
\end{multline}
\eqref{invega2}, \eqref{lkmleq} et \eqref{ts} donnent
$1/a-1/a_k \sim \frac{1}{(\ro-1)(k\log k)^{(\ro-1)}}$,
ce qui,  avec \eqref{invega}, termine la preuve.
\end{dem}

\begin{lem}\label{akpmap}
Il existe une constante positive $C_9$ telle que, pour tout $k \ge 2$,
et tout nombre premier $p$,
\begin{equation}\label{six13}
\abs{\frac{a_k}{p^{\ro_k}-1}  - \frac{a}{p^{\ro}-1}} 
\le C_9\frac{\log p}{p^{\ro_2}}\frac{1}{(k\log k)^{\ro-1}}\cdot
\end{equation}
\end{lem}

\begin{dem}
Par le théorème des accroissements finis appliqué à la
fonction $t \mapsto 1/(p^{t}-1)$, il existe 
$\theta,\ $  $\ro_k < \theta < \ro$,
tel que
\begin{eqnarray}
\abs{\frac{a_k}{p^{\rok}-1}-\frac{a}{p^{\ro}-1}}
&=& \left(\frac{1}{p^{\rok}-1}-\frac{1}{p^{\ro}-1} \right) a_k
+ \frac{a_k-a}{p^{\ro}-1}\nonumber\\
&=& \frac{p^{\theta}\log p}{(p^{\theta}-1)^2}(\ro-\rok)a_k
   + \frac{a_k-a}{p^\ro-1}\cdot
\label{akpa}
\end{eqnarray}
On a 
$1 < \ro_2 \le \ro_k \le \theta \le \ro \le 2$.
La décroissance de $x \mapsto x/(x-1)^2$ donne
\[
\frac{p^{\theta}}{(p^{\theta}-1)^2} \le 
\frac{p^{\ro_2}}{(p^{\ro_2}-1)^2} \le \frac{C}{p^{\ro_2}}
\]
où $C$ ne dépend ni de $p$ ni de $k$.
De même il existe  $D$ tel que
\[
\frac{1}{p^{\ro}-1} 
\le \frac{1}{p^{\ro_2}-1} 
\le \frac{D}{p^{\ro_2}}
\le \frac{D\log p}{p^{\ro_2}\log 2}
\le \frac{3 D\log p}{2 p^{\ro_2}}\cdot
\]
Et donc, puisque la suite $(a_k)$ est majorée par $3/2$,
par \eqref{akpa} on a

\[
\abs{\frac{a_k}{p^{\rok}-1}-\frac{a}{p^{\ro}-1}} \le
\frac{3\log p}{2p^{\ro_2}}(C(\ro-\rok) + D(a_k -a))
\]
et l'on conclut en utilisant les propositions
\ref{romrokprop} et \ref{akma}.
\end{dem}

\section{Un problème d'optimisation}

\subsection{La fonctions $F$}

\begin{definition}\label{defF}
Pour $\brx = (x_1,x_2,\dots,x_k) \in \Rplus^{*k}$ on définit
\begin{equation}\label{defFk}
F(\brx) = \sum_{j=1}^{k} x_j \log\left(1+\frac{c(\brx)}{x_j}\right) =
\sum_{j=1}^{k} (x_j \log(x_j + c(\brx)) - x_j \log x_j).
\end{equation}
\end{definition}

\begin{rem}
La fonction $F$ se prolonge par continuité sur $\R_+^k$ en
posant $0\log 0 = 0$. Notons que, par \eqref{cxk0}, on a
$F(x_1,x_2,\dots,x_k,0) = F(x_1,x_2,\dots,x_k)$ et
$F(\underline 0) = 0$.
\end{rem}

\begin{lem}
La fonction $F$ est concave dans $\Rplus^{k}$.
\end{lem}

\begin{dem}
Par \eqref{defFk} et \eqref{dronc} pour tout $x \in \Rplus^{*k}$,
on a
\begin{eqnarray}
\dron{F}{x_i}(\brx)
 &=& 
\left(\sum_{j=1}^{k}\frac{x_j}{\cxj} \right)\dron{c}{x_i}(\brx)
+  \log(\cxi) + \frac{x_i}{\cxi}-\log x_i-1\nonumber\\
&=& \frac{c(\brx)}{\cxi} + \frac{x_i}{\cxi} 
 + \log\left(\frac{\cxi}{x_i} \right) - 1\nonumber\\
&=& \log\left(\frac{\cxi}{x_i} \right).\label{d1F}
\end{eqnarray}
On a ensuite
\begin{eqnarray}
\ddron{F}{x_i}(\brx) &=& \dr{x_i}(\log(\cxi)-\log(x_i))\nonumber\\
&=& \left(\frac{1}{\cxi}\right)
\left(\dron{c}{x_i}(\brx) + 1\right) -\frac{1}{x_i}\nonumber\\
&=& \frac{c(\brx)}{(\cxi)^2\, T(x)} + \frac{1}{\cxi} - \frac{1}{x_i}
\end{eqnarray}
et, pour $j \ne i$,
\begin{eqnarray}
\dddron{F}{x_i}{x_j} &=& \dr{x_j}(\log(\cxi)-\log(x_i))\nonumber\\
&=&  \left(\frac{1}{\cxi}\right)\dron{c}{x_j}(\brx) 
= \frac{c(\brx)}{(\cxi)(\cxj)T(x)}\cdot
\end{eqnarray}
La forme quadratique des dérivées secondes de $F$ s'écrit
donc,
pour $\brx\in\Rplus^{\star k}$,
\begin{equation}\label{fqF}
F\sg(\brx)\cdot(h_1,h_2,\dots,h_k)  = 
\frac{c(\brx)}{T(\brx)}\left(\sum_{i=1}^{k}\frac{h_i}{\cxi}\right)^2
- \sum_{i=1}^{k} \frac{c(\brx)h_i^2}{x_i(\cxi)}\cdot
\end{equation}
L'inégalité de Cauchy-Schwarz donne 
\begin{eqnarray*}
\left(\sum_{i=1}^{k}\frac{h_i}{\cxi}\right)^2 &=& 
\left(\sum_{i=1}^{k}
\sqrt{\frac{x_i}{\cxi}}\frac{h_i}{\sqrt{x_i(\cxi)}}\right)^2\\
&\le& T(\brx) \sum_{i=1}^{k} \frac{h_i^2}{x_i(\cxi)},
\end{eqnarray*}
ce qui, avec \eqref{fqF}, prouve la concavité de $F$ dans
l'adhérence $\Rplus^{k}$ de $\Rplus^{\star k}$.
\end{dem}

\subsection{Proximité de  $\mathbf{\AA(\brx)}$ 
et de  $\mathbf{\textbf{exp}(F(\brx))}$.}
Le lemme suivant montre que $F$ est une assez bonne
approximation de $\log \AA$ (cf. définition \ref{definitionalpha}).

\begin{lem}
Soit $\brx \in \Rplus^{*k}$ ; on a
\begin{equation}\label{proxfa}
\AA(\brx) = \frac{1}{2\sqrt 2} \exp(F(\brx))
\prod_{j=1}^{k} \frac{1}{s(x_i)},
\quad\text{ où }\quad
s(x_i) = \frac{\Gamma(x_i+1)}{x_i^{x_i}e^{-x_i}}
\end{equation}
est le terme correctif de la formule de Stirling
$\Gamma(x+1) = x^x e^{-x} s(x)$, de l'ordre de grandeur de $\sqrt{2\pi x}$.
Plus précisément on a l'encadrement
\begin{equation}\label{stirling}
\sqrt{2\pi x} \le s(x) \le e \sqrt{x},
\quad x \ge 1.
\end{equation}
\end{lem}

\begin{dem}
De la définition \ref{definitionalpha}, il suit
\[
\AA(\brx) = \frac{1}{2\sqrt 2} \exp(-\Omega(\brx))
\prod_{j=1}^{k} \frac{(c(\brx) + x_i)^{x_i}}{\Gamma(x_i+1)}
= \frac{1}{2\sqrt 2} 
\prod_{j=1}^{k} \frac{(c(\brx) + x_i)^{x_i}}{x_i^{x_i}}
\frac{x_i^{x_i}e^{-x_i}}{\Gamma(x_i+1)}\cdot
\]
La formule \eqref{stirling} se réduit à l'encadrement classique
de $\Gamma(x+1)$, $x^x e^{-x}\sqrt{2\pi x} \le \Gamma(x+1) \le
x^xe^{-x}e\sqrt{x}$.
\end{dem}

\noindent
De la définition de $s(j)$ résulte immédiatement le lemme
suivant.
\begin{lem}\label{lemstirling}
Pour tout $j \in \N$, on a
\begin{equation}\label{lemj1j}
\frac{s(j+1)}{s(j)} = e\left(\frac{j}{j+1}\right)^j,
\end{equation}
et lorsque $j$ tend vers l'infini
\[
\frac{s(j+1)}{s(j)}=1+\frac{1}{2j} + \bigo{\frac{1}{j^2}}\cdot
\]
\end{lem}

\subsection{Maximisation de $F$}

Pour $k\ge 2$, entier, $p_k$ est le $k\iem$ nombre premier. 
Soit $A > 0$ réel, et  $k \ge 2$ on considère le domaine 
$\cd(A) \subset \Rplus^{*k}$, défini par
\begin{equation}\label{defD}
x_1 \log 2 + x_2 \log 3 + \cdots + x_k \log p_k \le A.
\end{equation}
Par \eqref{d1F} la fonction $F$ définie par \eqref{defFk} est 
croissante par rapport à chaque variable, le problème
d'optimisation
\begin{equation}\label{opt1}
\begin{cases}
\brx \in \cd(A)\\
\max F(\brx).
\end{cases}
\end{equation}
a donc même solution que le problème
\begin{equation}\label{opt2}
\begin{cases}
x_1 \log 2 + x_2 \log 3 + \cdots + x_k \log p_k = A\\
\max F(x_1,x_2,\dots,x_k).
\end{cases}
\end{equation}
Ce problème a été résolu dans \cite{Evans}, lemme 6, par la méthode
des multiplicateurs de Lagrange. Le multiplicateur de Lagrange 
est la constante $\rok$ définie en \eqref{defro}.

\begin{lem}\label{maxF}
L'unique solution
$\brx\str = (x_1\str,x_2\str,\dots,x_k\str)$
du problème \eqref{opt1} satisfait
\begin{equation}\label{xcont}
x_1\str \log 2 + x_2\str\log 3 + \cdots + x_k\str\log p_k = A,
\end{equation}
\begin{equation}\label{lagr}
\dron{F}{x_i}({\brx}\str) = \rok\log p_i.
\end{equation}
\begin{equation}\label{defak}
c(\bar{x}\str) = a_k A 
\end{equation}
\begin{equation}\label{defxis}
x_i\str = \frac{a_k A}{p_i^{\rok}-1},
\quad i=1,2,\dots, k
\end{equation}
\begin{equation}\label{maxvalue}
F({\brx}\str) = \rok A.
\end{equation}
\end{lem}

\begin{lem}\label{majFalpha}
Soit $k \ge 2$, $\bar \al = (\al_1,\al_2,\dots,\al_k) \in \cd(A)$
défini par \eqref{defD}, 
$\brx\str$ défini par \eqref{defxis} 
et $F$ définie par \eqref{defFk}.
Alors on a
\begin{eqnarray}
F(\bar \al) &\le& F(\brx\str) -
\frac{1}{4A\log p_k}\left(\sum_{i=1}^{k-1}\abs{\al_i-x_i\str}\log p_i
\right)^2\label{majFx}\\
&\le& F(\brx\str) - \frac{1}{4A\log p_k}
\sum_{i=1}^{k-1}(\al_i-x_i\str)^2(\log p_i)^2.\label{majFx2}
\end{eqnarray}
\end{lem}

\begin{dem}
Définissons $\bar y = (y_1,y_2,\dots,y_k)$ par
\begin{equation}\label{defy}
y_1 = \al_1,\ y_2 = \al_2,\dots,\ y_{k-1} = \al_{k-1},\quad\text{ et }\
\sum_{i=1}^{k} y_i \log p_i = A.
\end{equation}
Comme $ \al\in\cd(A)$, on a $\al_k \le y_k$ et la croissance
de $F$ par rapport à chacune des variables (cf. \eqref{d1F})
 entraine 
\begin{equation}\label{majxy}
F(\bar \al) \le F(\bar y).
\end{equation}
Posons $h_i = y_i - x_i\str$. On a par \eqref{defy}, \eqref{defxis}
et \eqref{defaak},
\begin{equation}\label{hyp}
\sum_{i=1}^{k} h_i \log p_i = \sum_{i=1}^{k} y_i \log p_i
- \sum_{i=1}^{k} x_i\str\log p_i = A-A = 0.
\end{equation}
La formule de Taylor appliquée à $F$ entre les points
$\bar y$ et $\brx\str$ donne
\begin{equation}\label{tlxy}
F(\bar y) - F(\brx\str)= 
\sum_{i=1}^{k} h_i\dron{F}{x_i}(\brx\str)
+ \frac{1}{2} F\sg(\underline\xi) \cdot (\bar h),
\end{equation}
avec $\bar\xi = \theta\bar y + (1-\theta)\brx\str$
et $0 < \theta < 1$.
On a donc $\bar\xi = (\xi_1,\xi_2,\dots,\xi_k) \in \Rplus^{\star k}$.
Par \eqref{lagr} et \eqref{hyp} il vient
\begin{equation}\label{hdfxz}
\sum_{i=1}^{k} h_i\dron{F}{x_i}(\brx\str) = 
\rok \sum_{i=1}^{k} h_i\log p_i = 0.
\end{equation}
\eqref{hyp} et l'inégalité de Cauchy-Schwarz donnent
\begin{eqnarray}
\left(\sum_{i=1}^{k} \frac{h_i}{\cxii} \right)^2
&=&
\left(\sum_{i=1}^{k}\frac{h_i}{\cxii}
\left(1-\frac{(\cxii)\log p_i}{2c(\bar \xi)\log p_k} \right)
\right)^2\nonumber\\
&=&
\left(\sum_{i=1}^{k}
\sqrt{\frac{\xi_i}{\cxii}}
\frac{h_i}{\sqrt{\xi_i(\cxii)}}
\left(1-\frac{(\cxii)\log p_i}{2c(\bar \xi)\log p_k} \right)
\right)^2\nonumber\\
&\le& T(\bar\xi)\label{lem362}
\left(
\sum_{i=1}^{k} \frac{h_i^2}{\xi_i(\cxii)}
\left(1-\frac{(\cxii)\log p_i}{2c(\bar \xi)\log p_k} \right)^2
\right)\cdot
\end{eqnarray}
Par \eqref{defcenc}, pour tout $i$, $1 \le i \le k$, on a
$0 < \xi_i \le \Omega(\bar \xi) \le c(\bar \xi)$.
En notant $t_i = \frac{(\cxii)\log p_i}{2c(\bar \xi)\log p_k}$,
on a donc $0 < t_i \le 1$, puis $(1-t_i)^2 \le 1-t_i$.
La majoration \eqref{lem362} et\\[2mm]
\eqref{fqF} donnent alors
\begin{eqnarray*}
F\sg(\bar\xi)\cdot \bar h &\le& 
\sum_{i=1}^{k}\frac{h_i^2 c(\bar\xi)}{\xi_i(\cxii)}
\left(1- \frac{(\cxii)\log p_i}{2c(\bar \xi)\log p_k} -1\right) \\
&\le& -\sum_{i=1}^{k} \frac{\log p_i}{\log p_k} \frac{h_i^2}{2\xi_i}\cdot
\end{eqnarray*}
Avec l'inégalité de Cauchy-Schwarz sous la forme
\[
\left(\sum_{i=1}^{k}\abs{h_i}\log p_i \right)^2 =
\left(\sum_{i=1}^{k} \sqrt{\xi_i\log p_i}
\frac{\abs{h_i}\sqrt{\log p_i}}{\sqrt\xi_i}\right)^2
\le
A \sum_{i=1}^{k} \frac{h_i^2\log p_i}{\xi_i}
\]
(car par \eqref{defy} et \eqref{xcont}, $\sum_{i=1}^{k}\xi_i\log p_i =
A$), cela donne
\[
F\sg(\xi)\cdot \bar h \le
-\frac{1}{2A\log p_k}\left(\sum_{i=1}^{k}\abs{h_i}\log p_i \right)^2
\le
-\frac{1}{2A\log p_k} \left(\sum_{i=1}^{k-1}\abs{h_i}\log p_i \right)^2
\]
ce qui avec \eqref{majxy}, \eqref{tlxy}, \eqref{hdfxz} et \eqref{defy}, 
complète la preuve du lemme \ref{majFalpha}.
\end{dem}

\begin{theorem}\label{thKnFk}
Il existe deux constantes $C_3$ et $C_4$ telles que, pour tout entier
$n \ge 2$, avec $n = 2^{\al_1}3^{\al_2}\dots p_k^{\al_k}$ avec
$\al_1,\al_2,\dots,\al_k \ge 1$ et $\bra = (\al_1,\al_2,\dots,\al_k)$
on ait
\begin{equation}\label{approxkf}
C_3 \frac{\exp(F(\bra))}{e^k\sqrt{\al_1 \al_2 \dots \al_{k}}}  \le K(n)
\le C_4 \frac{\exp(F(\bra))}{\pi^{k/2}}
\le C_4 \frac{\exp(\rok\log n)}{\pi^{k/2}}\cdot
\end{equation}
\end{theorem}

\begin{dem}
On part de l'encadrement donné par \eqref{ckc},
\[
C_1 \sqrt\pi\, \AA(\bar \al) \BB(\bar \al) 
\le K(n) \le 
 C_2\sqrt\pi\, \AA(\bar \al)\BB(\bar \al).
\]
Avec les encadrements de $\BB$, \eqref{encb}, et de $c$,
\eqref{defcenc}, ceci donne
\[
C_1 \sqrt{2\pi\Omega(\bra)}\AA(\bra)
\le K(n)
\le
C_{2}\sqrt{6\pi}\sqrt{\Omega(\bra)}\AA(\bra). 
\] 
On remplace $\AA(\bra)$ par le second membre de \eqref{proxfa},
\begin{multline*}
\frac{C_1}{2} \sqrt{\pi\Omega(\bra)}\exp(F(\bra)
 \prod_{j=1}^{k} \frac{1}{s(\al_j)}
\le K(n)\\ \le
\frac{C_{2}\sqrt{3\pi}}{2}
\sqrt{\Omega(\bra)}\exp(F(\bra) \prod_{j=1}^{k} \frac{1}{s(\al_j)}\cdot
\end{multline*}
L'encadrement \eqref{stirling} de $s(x)$ donne alors
\begin{multline}\label{l48bis}
\frac{C_1}{2} \sqrt{\pi\Omega(\bra)}\frac{1}{e^k \sqrt{\al_1\al_2\dots\al_k}.}
\exp(F(\bra))
\le K(n) \\
\le
\frac{C_{2}\sqrt{3\pi}}{2}
\sqrt{\frac{\Omega(\bra)}{\prod_{i=1}^{k}2\pi \al_i}}
\exp(F(\bra))\cdot
\end{multline}
La première inégalité dans
\eqref{approxkf} est obtenue en minorant $\Omega(\bra)$
par $1$ et en posant $C_3 = \frac{C_1\sqrt{\pi}}{2}$. 
Puisque chacun des entiers $2\al_i$ est supérieur ou égal à $2$,
leur somme $2\Omega(n)$ est majorée par leur produit
$\prod_{i=1}^{k} (2\al_i)$. On obtient donc la deuxième inégalité
de \eqref{approxkf} en choisissant $C_4 = \frac{C_2\sqrt{3\pi}}{2}$.

La dernière inégalité dans \eqref{approxkf} se réduit à $F(\bra) \le
\rok\log n$ ; elle s'obtient en appliquant le lemme \ref{maxF},
avec $A = \log n$, ce qui assure l'appartenance de $\bra$ au
domaine $\cd(A)$.
\end{dem}

\section{Grandes valeurs de la fonction $K$}

\begin{theorem}\label{grandval}
Il existe deux constantes
positives $C_5$ et $C_6$ telles que~:
\begin{enumerate}
\item
Pour tout entier $n$ suffisamment grand on a
\begin{equation}\label{grandval.1}
\log K(n)
\le \ro\log n - C_5 \frac{(\log n)^{1/\ro}}{\log\log n}\cdot
\end{equation}
\item
Pour tout $n$ suffisamment grand il existe $m \le n$ tel que
\begin{equation}\label{grandval.2}
\log K(m) 
\ge \ro\log n - C_6 \frac{(\log n)^{1/\ro}}{\log\log n}
\ge \ro\log m - C_6 \frac{(\log m)^{1/\ro}}{\log\log m}\cdot
\end{equation}
\end{enumerate}
\end{theorem}

\paragraph{Preuve de \eqref{grandval.1}}
Soit $n = q_1^{\al_1} q_2^{\al_2} \dots q_k^{\al_k}$ la décomposition de 
$n$ en facteurs premiers. 
Lorsque $k=1$, on calcule par récurrence avec la formule
\eqref{Kconv},
$K(n) = K(q_1^{\al_1}) = 2^{\al_1-1} \le n$
et \eqref{grandval.1} est vérifiée pour tout $C_5$ et $n$ assez grand.

Nous supposerons maintenant $k \ge 2$.
On pose  $N = p_1^{\al_1} p_2^{\al_2} \dots p_k^{\al_k} \le n$.
D'après la définition de $K$  on a $K(n) = K(N)$, et 
$\al_1\log 2 + \al_2 \log 3 + \cdots + \al_k \log p_k = \log N$.
L'encadrement \eqref{approxkf} donne
\[
\log K(n) = \log K(N) \le \rok\log N  - \frac{k}{2}\log \pi + \log C_4,
\]
qui, avec $\log N \le \log n$, donne
\begin{eqnarray}\nonumber
\log K(n) &\le& \rok\log n  - \frac{k}{2}\log \pi + \log C_4\\
&=& \ro\log n 
 - \left[(\ro -\rok)\log n +  \frac{k}{2}\log \pi -\log C_4\right].
\label{majlK}
\end{eqnarray}
Supposons $n \ge 16$ ce qui assure $\log\log n > 1$.
Vue \eqref{romrok}, il existe une constante positive
$\gamma_1$ telle que 
$\ro-\rok \ge \gamma_1/(k^{\ro-1}(\log k)^{\ro})$.
Alors
\begin{enumerate}
\item
Si \dsm{2 \le k \le \frac{(\log n)^{1/\ro}}{\log\log n} < \log n},
on a
\[
\ro-\rok 
\ge \frac{\gamma_1}{k^{\ro-1}(\log k)^{\ro}}
\ge 
\frac{\gamma_1}{\left(\frac{(\log n)^{1/\ro}}{\log\log n}\right)^{\ro-1}
(\log\log n)^{\ro}}
= \gamma_1 \frac{\left(\log n \right)^{1/\ro-1}}{\log\log n}\cdot
\]
\item
Si \dsm{k > \frac{(\log n)^{1/\ro}}{\log\log n}}, on a alors
\[
\frac{k}{2}\log \pi \ge \frac{k}{2} 
> \frac{1}{2}\frac{(\log n)^{1/\ro}}{\log\log n}\cdot
\]
\end{enumerate}
Dans les deux cas, le crochet de \eqref{majlK} vérifie
\begin{eqnarray*}
 \left[(\ro -\rok)\log n +  \frac{k}{2}\log \pi -\log C_4\right] 
&\ge& 
\min\left(\gamma_1,\frac{1}{2}\right)\frac{(\log n)^{1/\ro}}{\log\log n} -
\log C_4\\
&\ge& C_5 \frac{(\log n)^{1/\ro}}{\log\log n}
\end{eqnarray*}
pour $n$ assez grand, avec $C_5 > 0$.
Ceci termine la démonstration de \eqref{grandval.1}. 
Avant de démontrer \eqref{grandval.2}, rappelons le lemme suivant
qu'on trouvera dans \cite{Her-Nic}.
\begin{lem}\label{d2}
Soit $k$ un entier positif ; on range les $2^k$ diviseurs de
$n = p_1p_2\dots p_k$ par ordre croissant~: 
$1 = d_1 < d_2 < \dots < d_{2^k} = n$. Alors, pour tout 
$i$, $1 \le i \le 2^k-1$, on a $d_{i+1} \le 2 d_i$.
\end{lem}

\paragraph{Preuve de \eqref{grandval.2}}

On applique le lemme \ref{maxF} avec $A = \log n$ et 
\begin{equation}\label{defkapa}
k = \floor{\kappa \frac{(\log n)^{1/\ro}}{\log\log n}}
\end{equation}
où $\kappa$ est une constante positive satisfaisant
\begin{equation}\label{kpmaj}
\kappa < \ro a^{1/\ro} = 1.82\dots
\end{equation}
et $a$ est défini en \eqref{defaak}.
Par le lemme \ref{maxF} le maximum de $F$ est atteint
en 
\begin{equation}\label{xistar}
{\bar x}\str = (x_1\str,x_2\str,\dots,x_k\str)
\quad\text{ avec }\quad
x_i^{\star} = \frac{a_k \log n}{p_i^{\rok}-1},\quad i=1,2,\dots,k,
\end{equation}
où les $a_k$ sont donnés par \eqref{defaak}. On rappelle que,
par \eqref{xcont},
\begin{equation}\label{somxlogp}
\sum_{i=1}^{k} x_i^{\star} \log p_i = \log n.
\end{equation}
Lorsque $n \to \infty$, on a, en utilisant la décroissance
de $(a_k)$ (cf. lemme \ref{akma}),
puis \eqref{defkapa} pour
obtenir un équivalent de $k$,
\begin{equation}\label{seqxistar}
x_1^{\star} > x_2^{\star} > \cdots > x_k^{\star} = 
\frac{a_k \log n}{p_k^{\rok}-1} > \frac{a\log n}{p_k^{\ro}}
\sim\frac{a\log n}{(k\log k)^{\ro}}
\sim \frac{a\ro^\ro}{\kappa^{\ro}} > 1,
\end{equation}
la dernière inégalité provenant de \eqref{kpmaj}.

\paragraph{Construction de $m$.}

$k$ étant défini par \eqref{defkapa} et $\bar x$ par \eqref{xistar} on
pose 
$m_0 = \prod_{i=1}^{k} p_i^{\floor{x_i^{\star}}}$.
Par \eqref{somxlogp} on a
\[
\frac{n}{p_1p_2\dots p_k} = \prod_{i=1}^{k} p_i^{x_i^{\star}-1}
< m_0 \le \prod_{i=1}^{k} p_i^{x_i^{\star}} = n.
\]
Soit $d$ le plus grand diviseur de $p_1p_2\dots p_k$
vérifiant $d \le n/m_0 < p_1 p_2 \dots p_k$. 
On pose $m = m_0 d$. 
Par définition de $d$ et par le lemme \ref{d2} on a 
$n/m_0 < 2d$, et donc
\begin{equation}\label{encadm}
1 \le\frac{n}{m} = \frac{n}{m_0 d} <  2.
\end{equation}
On écrit  $d = \prod_{i=1}^{k} p_i^{\veps_i}$, 
avec $\veps_i \in \set{0,1}$. 
On a alors
\begin{equation}\label{defm}
m = \prod_{i=1}^{k} p_i^{\al_i}
\end{equation}
avec 
\begin{equation}\label{alphai}
\al_i = \floor{x_i^{\star}} + \veps_i,\quad \veps_i \in \set{0,1}.
\end{equation}
Avec \eqref{seqxistar}, \eqref{alphai} donne
\begin{equation}\label{minxi}
1 \le \floor{x_i\str} \le \al_i \le \floor{x_i\str} + 1
\le x_i\str + 1 \le 2 x_i\str, \quad 1 \le i \le k.
\end{equation}
et, par \eqref{encadm}, \eqref{somxlogp} et \eqref{defm} il vient
\begin{equation}\label{minlogmn}
-\log 2 \le \log\frac{m}{n} = \sum_{i=1}^{k} (\al_i - x_i\str)\log p_i
\le 0.
\end{equation}

\paragraph{Fin de la preuve de \eqref{grandval.2}}
La formule de Taylor et l'expres\-sion \eqref{fqF} de la forme
quadratique dérivée seconde de $F$ donnent
\begin{multline}\label{tl2}
F(\bar\al) = F(\bar{x^{\star}}) +
 \sum_{i=1}^{k} (\al_i - x_i^{\star})\dron{F}{x_i}(\bar{x^{\star}})
+ \frac{c(\bar\xi)}{T(\bar\xi)}
   \left(\sum_{i=1}^{k}\frac{\al_i-x_i^{\star}}{\cxii}\right)^2\\
- \sum_{i=1}^{k}
  \frac{c(\bar\xi)(\al_i - x_i^{\star})^2}{\xi_i(c(\bar\xi)+\xi_i)}
\end{multline}
avec, pour $1 \le i \le k$,\ 
$\xi_i = \theta\al_i + (1-\theta)x_i^{\star}$,
$0 < \theta < 1$ et $\bar\xi = (\xi_1,\xi_2,\dots,\xi_k)$.
Par \eqref{minxi} on a
\begin{equation}
\xi_i \ge \theta \floor{x_{i}^{\star}}  + (1-\theta)\floor{x_{i}^{\star}} 
= \floor{x_{i}^{\star}} \ge 1. 
\end{equation}
Par \eqref{alphai}, le quatrième terme du second membre de \eqref{tl2}
satisfait donc
\begin{equation}\label{tm4}
- \sum_{i=1}^{k} 
\frac{c(\bar\xi)(\al_i - x_i^{\star})^2}{\xi_i(c(\bar\xi)+\xi_i)}
\ge - \sum_{i=1}^{k} (\al_i - x_i^{\star})^2 \ge -k.
\end{equation}
Le troisième terme est positif.  Le second terme, par \eqref{lagr},
\eqref{minlogmn} puis le lemme \ref{derivc}
satisfait
\begin{multline}\label{tm2}
\sum_{i=1}^{k} (\al_i - x_i^{\star})\dron{F}{x_i}(\bar x^{\star})
= \rok \sum_{i=1}^{k} (\al_i - x_i^{\star}) \log p_i
= \rok \log\frac{m}{n}\\
\ge -\rok\log 2 \ge -\ro\log 2.
\end{multline}
De \eqref{tl2}, \eqref{tm4} et \eqref{tm2} on déduit
\begin{equation}\label{ttt}
F(\bar \al) \ge F(\bar x\str) - \ro\log 2 -k.
\end{equation}
Vue la décomposition en facteurs premiers de $m$
\eqref{defm},  l'encadrement \eqref{approxkf} donne
\begin{equation}
\log K(m) \ge F(\bar \al) 
 - k -\frac12\sum_{i=1}^{k} \log \al_i + \log C_3.
\end{equation}
Avec \eqref{ttt} cela donne
\begin{equation}\label{}
\log K(m) \ge F(\bar x\str)
-2k -\frac12\sum_{i=1}^{k} \log \al_i + \log C_3
-\ro\log 2.
\end{equation}
Par \eqref{maxvalue} il vient
\begin{equation}\label{fxbar}
F(x_1^{\star},x_2^{\star},\dots,x_k^{\star}) 
= \ro\log n - (\ro-\ro_k)\log n,
\end{equation}
d'où
\begin{equation}\label{toto}
\log K(m) \ge  \ro\log n - (\ro-\ro_k)\log n -2k
-\frac12\sum_{i=1}^{k} \log \al_i + \log C_3 -\ro\log 2.
\end{equation}
Il reste à montrer que, à part  $\ro\log n$,
les termes du second membre de \eqref{toto} sont tous des
$O(\log(n)^{1/\ro}/\log\log n) = O(k)$.
Pour le terme $(\ro-\rok)\log n$ cela résulte de l'équivalence
\eqref{romrok} et de \eqref{defkapa} 
Il reste donc à démontrer que
$\sum_{i=1}^{k} \log\al_i$ est un $O(k)$ c'est-à dire que
$\frac{1}{k}\sum_{i=1}^{k} \log\al_i$ est un $O(1)$
lorsque $k$ tend vers l'infini.
Par \eqref{minxi},
\[
\frac{1}{k}\sum_{i=1}^{k} \log \al_i 
\le \frac{1}{k} \sum_{i=1}^{k}\log(2x_i\str)
\le \log 2 + \frac{1}{k}\sum_{i=1}^{k} \log x_i\str.
\]
Vu le lemme suivant, cela termine la preuve de
\eqref{grandval.2} et du théorème \ref{grandval}.
\hfill$\Box$

\begin{lem}
Avec le choix de $k$, \eqref{defkapa}, on a, lorsque $n$ tend vers l'infini
\[
\frac{1}{k}\sum_{i=1}^{k} \log x_i\str = \bigo{1}.
\]
\end{lem}

\begin{dem}
Puisque $A = \log n$, la formule \eqref{xistar} donne
\begin{eqnarray*}
\frac{1}{k}\sum_{i=1}^{k} \log x_i\str &=& \log a_k + \log\log n
- \frac{1}{k} \sum_{i=1}^{k}\log(p_i^{\ro_k}-1)\\
&=&  \log a_k + \log\log n
- \sum_{i=1}^{k}\frac{\ro_k}{k}\log p_i 
- \frac{1}{k}\sum_{i=1}^{k}\log\left(1-\frac{1}{p_i^{\ro_k}}\right)\\
&=&  \log\log n - \frac{\ro_k}{k} \sum_{i=1}^{k}\log p_i + O(1).
\end{eqnarray*}

En utilisant le développement asymptotique de $p_k$
(cf. \cite{landau}, \S 57)
\[
p_k = k(\log k + \log\log k + \bigo{1})
\]
et le théorème des nombres premiers sous la forme
$\sum_{p \le x} \log p = x + \bigo{x/\log x}$
(cf. \cite{emf}, théorème 4.7) on obtient
\[
\sum_{i=1}^{k} \log p_i = p_k + \bigo{\frac{p_k}{\log p_k}}
=  k(\log k + \log\log k + \bigo{1}).
\]
D'où
\begin{equation}\label{ouf}
\frac{1}{k}\sum_{i=1}^{k} \log x_i\str = 
\log\log n - \ro_k[\log k + \log\log k]
+ O(1).
\end{equation}
En notant $\log_3 x = \log\log\log x$,
il résulte de la définition \eqref{defkapa} de $k$ que
\[
\log k = \frac{1}{\ro}\log\log n - \log_3 n + O(1),
\qquad
\log\log k = \log_3 n + O(1),
\]
et $\ro(\log k + \log\log k) = \log_2 n + \bigo{1}$.
Avec \eqref{ouf} cela donne
\[
\frac{1}{k}\sum_{i=1}^{k} \log x_i\str = 
(\ro-\rok)(\log k + \log\log k) + O(1) = O(1)
\]
car, par l'équivalence  \eqref{romrok}, 
$(\ro-\rok)(\log k+\log\log k)$
tend vers $0$ quand $n$ tend vers l'infini.
\end{dem}

\section{Propriétés des nombres $K$-champions}

Comme la fonction $K$ ne dépend que des exposants de
la décomposition en facteurs premiers de $n$, il est
clair que tout nombre $K$-champion $N$ est de la forme
\begin{equation}\label{factorchamp}
N = 2^{\al_1}3^{\al_2}\dots p_k^{\al_k},
\quad\text{ avec }\quad
\al_1 \ge \al_2 \ge \cdots \ge \al_k \ge 1.
\end{equation}
Il résulte de l'inégalité \eqref{grandval.2} du théorème \ref{grandval} que
\[
\limsup_{n \to \infty} K(n) = +\infty,
\]
et donc il existe une infinité de nombres $K$-champions.
$K(n)$ est le nombre de solutions en entiers
plus grands que $1$ de
\[
n = d_1 d_2 \dots d_r.
\]
Chaque telle factorisation de $n$ donne une factorisation
non triviale de $2n$, $2n = 2 d_1 d_2 \dots d_r$. Comme $2n$ admet
aussi la factorisation triviale $2n = (2n)$,
on a, pour $n \ge 2$, $K(2n) > K(n)$.
Il en résulte
que si l'on range les nombre $K$-champions dans
l'ordre  croissant
\begin{equation}\label{septantedeuxbis}
N_1 = 1 < N_2 = 4 < N_3 = 6 < \cdots < N_i < N_{i+1} < \cdots
\end{equation}
on a 
\begin{equation}\label{facteur2}
N_{i+1} \le 2 N_i.
\end{equation}

\subsection{Encadrement de $\omega(N)$}

\begin{theorem}\label{thomega}
\mbox{}

\begin{enumerate}
\item
Soit $N$ un nombre $K$-champion assez grand. Alors
\begin{equation}\label{thomega1}
\log K(N) \ge \ro\log N - C_6 \frac{(\log N)^{1/\ro}}{\log\log N}
\end{equation}
où $\ro$ est défini en \eqref{defro} et $C_6$ dans
le théorème \ref{grandval}.
\item
De plus, il existe trois constantes positives
$C_7$, $C_8$ et $N_0$ telles
que, pour tout $K$-champion $N \ge N_0$, on ait
\begin{equation}\label{thomega2}
C_7 \frac{(\log N)^{1/\ro}}{\log\log N}
\le \omega(N)
\le  C_8 \frac{(\log N)^{1/\ro}}{\log\log N}\cdot
\end{equation}
\end{enumerate}
\end{theorem}

\paragraph{Preuve de \eqref{thomega1}}
Pour le premier point, on utilise l'inégalité
\eqref{grandval.2} du théorème \ref{grandval}.
Il existe $m \le N$ vérifiant
$K(m) \ge \ro\log N - C_6 (\log N)^{1/\ro}/\log\log N$. 
Puisque $N$ est un $K$-champion, on a $K(N) \ge K(m)$.

\paragraph{Preuve de \eqref{thomega2}}
Posons $k = \omega(N)$.
L'inégalité \eqref{approxkf} donne
\begin{eqnarray}
\log K(N) &\le& \rok\log N - k \frac{\log \pi}{2} + \log C_4\nonumber\\
&\le& \ro\log N - (\ro-\rok)\log N
 - k \frac{\log \pi}{2} + \log C_4 .\label{tutu}
\end{eqnarray}
Vue \eqref{thomega1} cela donne, quand $N\to+\infty$,
\begin{equation}
(\ro-\rok)\log N  \le C_6 \frac{(\log N)^{1/\ro}}{\log\log N}
+\log C_4
\le  C_6(1+o(1)) \frac{(\log N)^{1/\ro}}{\log\log N}\cdot
\end{equation}
Par l'équivalence \eqref{romrok}, il existe donc $C > 0$ tel que,
\begin{equation}\label{th41}
k^{\ro-1}(\log k)^{\ro} \ge C (\log N)^{1-1/\ro}\log\log N.
\end{equation}
Mais la fonction $y = f(t) = t^{\ro-1}(\log t)^{\ro}$ 
tend vers l'infini avec $t$, et est croissante pour $ t\ge 1$.
Sa fonction réciproque satisfait
\[
f^{-1}(y) \sim (\ro-1)^{\frac{\ro}{\ro-1}}
\left(\frac{y}{(\log y)^{\ro}}\right)^{\frac{1}{\ro-1}}
\asymp \left(\frac{y}{(\log y)^{\ro}}\right)^{\frac{1}{\ro-1}},
\qquad (y \to \infty).
\]
et \eqref{th41} entraîne
\[
k \ge f^{-1}\left(C(\log N)^{1-\frac{1}{\rho}}\log\log N \right)
\asymp \frac{(\log N)^{1/\rho}}{\log\log N}\cdot
\]
Ceci donne la minoration de $k=\omega(N)$ dans \eqref{thomega2}.
La comparaison de \eqref{tutu} et \eqref{thomega1}
donne aussi
\[
k \frac{\log\pi}{2} \le C_6 \frac{(\log N)^{1/\ro}}{\log\log N} + \log
C_4 \ll  \frac{(\log N)^{1/\ro}}{\log\log N}
\]
d'où l'on déduit la majoration de $\omega(N)$ dans \eqref{thomega2}.
\hfill$\Box$

\begin{lem}\label{somlog}
  Soit $N$ un nombre $K$-champion dont la décomposition en facteurs
  premiers est donnée par \eqref{factorchamp}.  On applique le
  lemme \ref{maxF} avec $k = \omega(N)$ et $A = \log N$.  Soit
  $\bar x\str$ défini par \eqref{defxis}.  Alors, lorsque $N \to
  +\infty$, on a
\begin{equation}\label{lpk}
\log p_k \sim \frac{1}{\ro} \log\log N,
\quad\text{où $\ro$ est défini en \eqref{defro}},
\end{equation}
et
\begin{equation}\label{defdelta}
\sum_{i=1}^{k} \abs{\alpha_i -x_i\st}\log p_i = \bigo{(\log N)^{\delta}}
\quad\text{avec}\quad
\delta = (1+1/\ro)/2= 0.789\, 243\dots 
\end{equation}
\end{lem}

\begin{dem}
  L'inégalité \eqref{thomega1} du théorème \ref{thomega} et $\ro
  \ge \rok$ (cf. lemme \ref{romrokprop}) donnent
\begin{equation}\label{majlkN}
\log K(N) \ge \ro \log N  - C_6 \frac{(\log N)^{1/\ro}}{\log\log N}
\ge \rok \log N  - C_6 \frac{(\log N)^{1/\ro}}{\log\log N}\cdot
\end{equation}
L'équation \eqref{maxvalue} du lemme \ref{maxF}
donne $F(\bar x\str) = \rok\log N$.
Appliquons maintenant le  lemme \ref{majFalpha},
avec $\bar\al = (\al_1,\al_2,\dots,\al_k)$.
La majoration \eqref{majFx} s'écrit
\begin{equation}\label{minFaN}
F(\bar\al) \le \rok\log N -
\frac{1}{4\log\!N\,\log p_k}
\left(\sum_{i=1}^{k-1}\abs{\al_i-x_i\str}\log p_i\right)^2\cdot
\end{equation}
Par l'encadrement \eqref{thomega2} du théorème \ref{thomega}, $k =
\omega(N)$ tend vers l'infini avec $N$.  En utilisant
\eqref{approxkf}, il en résulte que, pour $N$ assez grand,
\[
\log K(N) \le  F(\bar\al) + \log C_4 - \frac{k}{2}\log \pi
\le F(\bar\al).
\]
Cette inégalité, avec \eqref{majlkN} et \eqref{minFaN} donne
\begin{equation}\label{sixneuf}
\frac{1}{4\log\!N\,\log p_k}
\left(\sum_{i=1}^{k-1}\abs{\al_i-x_i\str}\log p_i\right)^2
\le
C_6 \frac{(\log N)^{1/\ro}}{\log\log N}\cdot
\end{equation}
Par les inégalités \eqref{thomega2} du théorème \ref{thomega},
$k = \omega(N) \asymp \frac{(\log N)^{1/\ro}}{\log\log N}$,
ce qui entraîne, lorsque $N \to +\infty$,
\[
\log p_k \sim \log(k \log k) \sim \log k \sim \frac{1}{\ro}\log\log N.
\]
et cela prouve \eqref{lpk}.
De \eqref{sixneuf} et \eqref{lpk} on déduit
\begin{equation}\label{somlogpkm1}
\sum_{i=1}^{k-1} \abs{\al_i - x_{i}\str}\log p_i = \bigo{(\log N)^{\delta}}
\end{equation}
avec $\delta$  donné par \eqref{defdelta}. Pour prouver
que $\sum_{i=1}^{k}\abs{\al_i-x_{i}\str}\log p_i =\bigo{(\log
  N)^{\delta}}$,
ce qui terminera la preuve du lemme,
il reste à vérifier que
$\abs{(\al_k-x_k\str)\log p_k} = \bigo{(\log N)^{\delta}}$. 
La définition $A = \log N$ et \eqref{xcont} donnent
\[
\log N = \sum_{i=1}^{k} \al_i \log p_i = 
\sum_{i=1}^{k} x_i\str \log p_i.
\]
Il en résulte
\begin{equation}\label{somlogzero}
\sum_{i=1}^{k} (\al_i-x_i\str)\log p_i = 0
\end{equation}
et donc, avec \eqref{somlogpkm1}, 
\dsm{
\abs{(\al_k-x_k\str)\log p_k} \le 
\sum_{i=1}^{k-1} \abs{\al_i - x_{i}\str}\log p_i = \bigo{(\log N)^{\delta}}.
}
\end{dem}

\subsection{Exposants des petits facteurs premiers}

\begin{theorem}\label{decNth1}
Soit $N$ un nombre $K$-champion dont la décomposition
en facteurs premiers est donnée par \eqref{factorchamp},
et $\ro$ et $a$ définis en \eqref{defro} et \eqref{defaak}.
Lorsque $N \to +\infty$, on a
\begin{equation}\label{decNth1.2}
\Omega(N) = \al_1 + \al_2 + \cdots + \al_k = b\log N +
O\left((\log N)^\delta\right),
\end{equation}
avec $\delta$ défini en \eqref{defdelta},
\begin{equation}\label{defbetai}
\beta_i = \frac{a}{p_{i}^{\ro}-1}
\quad\text{ et }\quad
b = \sum_{i=1}^{+\infty} \beta_i = 0.8612985\dots
\end{equation}
De plus, pour $1 \le i \le k$, on a uniformément en $i$,
\begin{equation}\label{decNth1.1}
\alpha_i =  \beta_i \log N +
O\left(\frac{(\log N)^{\delta}}{\log p_i} \right).
\end{equation}
\end{theorem}

\paragraph{Démonstration de \eqref{decNth1.2}}
Pour prouver \eqref{decNth1.2} il suffit de démontrer la majoration
un peu plus forte
\begin{equation}\label{majfort}
\sum_{i=1}^{k}\abs{\al_i - \beta_i\log N}
+ \sum_{i=k+1}^{+\infty}\beta_i\log N = \bigo{(\log N)^\delta}.
\end{equation}
On applique le lemme \ref{maxF} avec $k=\omega(N)$
et $A=\log N$.
On écrit alors,
\[
\sum_{i=1}^{k} \abs{\al_i-\beta_i\log N} \le
\sum_{i=1}^{k} \abs{\al_i-x_i\st} +
\sum_{i=1}^{k} \abs{x_i\st-\beta_i\log N}
\]
où les $x_i\st$ sont définis par \eqref{defxis}.
Le lemme \ref{somlog} donne 
\begin{equation}\label{}
\sum_{i=1}^{k} \abs{\al_i-x_i\st} \le \bigo{(\log N)^\delta}.
\end{equation}
En utilisant le lemme \ref{akpmap} (car, par \eqref{thomega2}, $k \ge 2$)
\begin{eqnarray}\nonumber
\sum_{i=1}^{k} \abs{x_i\st-\beta_i \log N} &=&
\log N\sum_{i=1}^{k} \abs{\frac{a_k}{p_i^{\rok}-1}-\frac{a}{p_i^{\ro}-1}}\\
&\le& C_9 \frac{\log N}{(k\log k)^{\ro-1}} 
\sum_{i=1}^{k} \frac{\log p_i}{p_i^{\ro_2}}
= \bigo{\frac{\log N}{k^{\ro-1}}}\label{xistmbi}\cdot
\label{majglob}
\end{eqnarray}
Par la minoration \eqref{thomega2} de $\omega(N)=k$, il vient
\begin{equation}\label{majlognk}
\frac{\log N}{k^{\ro-1}} = \bigo{(\log N)^{1/\ro}(\log\log N)^{\ro-1}}
= \bigo{(\log N)^{\delta}}
\end{equation}
car $1/\ro < \delta$.
Il vient ensuite
\begin{equation}\label{cinq22bis}
\sum_{i=k+1}^{+\infty} \beta_i = 
\sum_{i=k+1}^{+\infty}\frac{a}{p_i^{\ro}-1}
\le \sum_{i=k+1}^{\infty}\frac{2a}{p_i^\ro} 
\le \sum_{i=k+1}^{\infty}\frac{2a}{i^\ro}
\le \int_{k}^{\infty} \frac{2adt}{t^\ro}
= \frac{2a}{(\ro-1)k^{\ro-1}}
\end{equation}
ce qui, avec \eqref{xistmbi} et \eqref{majlognk} complète
la preuve de \eqref{majfort} et \eqref{decNth1.2}.

\paragraph{Démonstration de \eqref{decNth1.1}}

Pour tout $i$, $1\le i \le k$, on a
\begin{equation}\label{six12}
\abs{\al_i - \beta_i\log N}  \le \abs{\al_i-x_i\st} +
\abs{x_i\st-\beta_i\log N}.
\end{equation}
Par le lemme \ref{somlog} on a
\begin{equation}\label{almxi}
\abs{\al_i-x_i\st} \ll \frac{1}{\log p_i}(\log N)^{\delta}. 
\end{equation}
Les définitions \eqref{defxis} et \eqref{defbetai} de $\beta_i$ et
$x_i\str$ et le 
lemme \ref{akpmap}  donnent, pour tout $i$,
\[
\abs{x_i\st-\beta_i\log N} \le 
C_9 \frac{\log p_i}{p_i^{\ro_2}}\frac{\log N}{(k\log k)^{\ro-1}}
= \frac{1}{\log p_i}\bigo{\frac{\log N}{k^{\ro-1}}}.
\]
En utilisant \eqref{majlognk}, il vient
\[
\abs{x_i\st-\beta_i\log N} \ll \frac{1}{\log p_i}(\log N)^{\delta},
\]
qui, avec \eqref{almxi} et \eqref{six12}, termine la preuve
de \eqref{decNth1.1}.
\hfill$\Box$

\subsection{Exposants des grands facteurs premiers}

\begin{lem}\label{climit}
Soit $N = 2^{\al_1} 3^{\al_2} \dots p_{k}^{\al_{k}}$ un $K$-champion
tendant vers l'infini,
$\bra = (\al_1,\al_2,\dots,\al_k)$, 
et $\delta$ défini par \eqref{defdelta}. 
\begin{enumerate}
\item
Soit $a = 1.100020011\dots$  la constante définie en \eqref{defaak};
alors
\begin{equation}\label{cbeta}
c(\bra) = a \log N   + \bigo{(\log N)^{\delta}},
\end{equation}
\item
Soit $T_0 = \sum_{j=1}^{\infty} \frac{1}{p_{j}^{\ro}} = 0.62035\dots$.
Alors on a
\begin{equation}\label{tbeta}
T(\bra) = T_0 +  \bigo{(\log N)^{\delta-1}}.
\end{equation}
\item
Soit $B_{0} = \sqrt{\frac{2a}{T_{0}}}= 1.883\dots$.
Alors on a
\begin{equation}\label{bbeta}
\BB(\bra) = B_{0} \sqrt{\log N}
\left(1+\bigo{(\log N)^{\delta-1}}\right).  
\end{equation}
\item
\begin{equation}\label{Kchampionequiv}
K(N) =  B_{0}\sqrt\pi \sqrt{\log N}\,\, \AA(\bra)
\left(1+\bigo{(\log N)^{\delta-1}} \right).
\end{equation}
\end{enumerate}
\end{lem}

\paragraph{Démonstration de \eqref{cbeta}}
Soit $\bar \beta = (\beta_i)_{i \ge 1}$ où les $\beta_i$ sont
définis en \eqref{defbetai}. L'égalité \eqref{chomogen}, 
le lemme \ref{lemcxy} et la majoration \eqref{majfort} donnent
\begin{eqnarray*}
\lefteqn{
\abs{c(\bra) - c(\bar\beta)\log N} = \abs{c(\bra)-c(\bar\beta \log N)}
\le}
\\
&\le&    2\sum_{i=1}^{k} \abs{\al_i - \beta_i\log N}
  + 2 \sum_{i=k+1}^{+\infty} \beta_i\log N =
\bigo{(\log N)^{\delta}}.
\end{eqnarray*}
Il reste à vérifier que $c(\bar\beta) = a$.
Or, par définition de $\ro$,
\[
2 = \prod_{i=1}^{\infty} \frac{1}{1-\frac{1}{p_{i}^{\ro}}} 
=  \prod_{i=1}^{\infty} \frac{p_i^\ro}{p_{i}^{\ro}-1} 
=  \prod_{i=1}^{\infty}\left(1+  \frac{1}{p_{i}^{\ro}-1}\right) \cdot
\]
Par la formule \eqref{defceq} qui définit $c$, cela s'écrit 
\dsm{
c\left(\frac{1}{2^{\ro}-1},\frac{1}{3^{\ro}-1},\cdots\right) = 1.
}
Puisque $\beta_i = a/(p_{i}^{\ro}-1)$, on obtient
$c(\bar\beta) = a$ en multipliant les deux termes de cette
égalité par $a$ (grâce à \eqref{chomogen}).

\paragraph{Démonstration de \eqref{tbeta}}
Par la défintion \eqref{defbetai}, on a
$\beta_i/(a+\beta_i) = 1/p_i^\ro$.
L'estimation \eqref{decNth1.1}, où le $\bigo{(\log N)^{\delta-1}}$
est uniforme en $i$, et \eqref{cbeta} donnent alors
\begin{eqnarray*}
T(\bra) &=& \sum_{i=1}^{k}\frac{\al_i}{c(\bra)+\al_i}
= \sum_{i=1}^{k}
\frac{\beta_i\log N \left(1+\bigo{(\log N)^{\delta-1}}\right)}
{(a+\beta_i)\log N\left(1+\bigo{(\log N)^{\delta-1}}\right)}
\\
&=& \sum_{i=1}^{k} 
\frac{1}{p_i^\ro}
\left(1+\bigo{(\log N)^{\delta-1}}\right)
= \sum_{i=1}^{k} \frac{1}{p_i^\ro} + \bigo{(\log N)^{\delta-1}}\\
&=& T_0 - \sum_{i=k+1}^{\infty} \frac{1}{p_{i}^{\ro}}
+ \bigo{(\log N)^{\delta-1}},
\end{eqnarray*}
ce qui prouve \eqref{tbeta}, car, par \eqref{majlognk} et
\eqref{cinq22bis}, on a
\[
\sum_{i=k+1}^{\infty} \frac{1}{p_i^{\ro}} = \bigo{\frac{1}{k^{\ro-1}}}
= \bigo{(\log N)^{\delta-1}}.
\]

\paragraph{Démonstration de \eqref{bbeta}}
La formule \eqref{bbeta} est une conséquence de la définition
\eqref{defB} de $\BB$, de \eqref{cbeta} et de \eqref{tbeta}.
\paragraph{Démonstration de \eqref{Kchampionequiv}}
Choisissant $\eta = 1-\delta$ dans le 
théorème \ref{evans} (théorème de Evans) on obtient
\begin{equation*}
K(N) = 
\sqrt{\pi}\BB(\bra)\AA(\bra)\left(1+\bigo{(\Omega(N))^{\delta-1}}\right)
\cdot
\end{equation*}
On conclut avec \eqref{bbeta} et \eqref{decNth1.2}.
\hfill$\Box$

\begin{definition}\label{defvoisin}
  Soit $N = 2^{\al_1}3^{\al_2}\dots p_k^{\al_k}$ un $K$-champion
  tendant vers l'infini, et $M$ un entier dépendant de $N$, $M =
  2^{\al_1'}3^{\al_2'}\dots p_{k'}^{\al_k'}$, avec $\al_i' \ge 1$
  pour $1 \le i \le k'$. On pose
  $\bra' = (\al_1',\al_2',\dots,\al'_{k'})$,
  et $\mu = \delta-1/\ro = 0.210\dots$.  
  On dira que $M$ est
  \emph{voisin} de $N$ si, avec la notation 1,
\[
\nsup{\bra'-\bra} = 
\bigo{(\log N)^\mu} = \bigo{(\log N)^{\delta-1/\ro}}.
\]
\end{definition}

\begin{lem}
Si $M = 2^{\al_1'}3^{\al_2'}\dots p_{k'}^{\al_k'}$ est un voisin de
$N$, on a
\begin{equation}\label{lembiz}
\log M = \log N + \bigo{(\log N)^{\delta}}.
\end{equation}
\begin{equation}\label{bigomegaequiv}
\Omega(M) = b\log N + O\left((\log N^\delta\right),\quad
\alpha_i' =  \beta_i \log N + O\left((\log N)^\delta\right)
\quad(1 \le i \le k').
\end{equation}
\end{lem}

\paragraph{Démonstration de \eqref{lembiz}}

Posons $k_1 = \max(k,k')$. On a
\begin{eqnarray*}
\abs{\log M-\log N} &=& \abs{\sum_{i=1}^{k_1} (\al_i'-\al_i)\log p_i}\\
&\le& \sum_{i=1}^{k_1} \abs{\al_i'-\al_i}\log p_i
\le \nsup{\bra'-\bra}\log p_{k_1}.
\end{eqnarray*}
Par \eqref{thomega2}, 
$k = \omega(N) = \bigo{\frac{(\log N)^{1/\ro}}{\log\log N}}$.
Puisque $M$ est voisin de $N$, on a $0 \le k_1-k \ll (\log N)^\mu$.
Par \eqref{thomega2}, on a donc $k_1 = \omega(N) + \bigo{(\log N)^\mu}
 = \bigo{\frac{(\log N)^{1/\ro}}{\log\log N}}$.
Par le théorème des nombres pre\-miers 
$p_{k_1} \sim k_1 \log k_1 = \bigo{(\log N)^{1/\ro}}$,
et, par hypothèse, $\nsup{\bra'-\bra} = \bigo{(\log
N)^{\delta-1/\ro}}$,
d'où \eqref{lembiz}.

\paragraph{Démonstration de \eqref{bigomegaequiv}}
L'égalité \eqref{decNth1.2} et la majoration évidente
\[
\abs{\Omega(M)-\Omega(N)} \le \nsup{\bra'-\bra}
= \bigo{(\log N)^{\delta}}
\]
donnent la première égalité dans \eqref{bigomegaequiv}. 
Pour $i$ satisfaisant  $i \le \min(k,k')$, l'estimation
\eqref{decNth1.1} et 
$\abs{\al_i'-\al_i}\, \le \nsup{\bra'-\bra} = \bigo{(\log
N)^{\delta}}$
donnent la deuxième égalité.
Lorsque $k < i \le k'$, on écrit
\[
\abs{\al_i'-\beta_i \log N} \le \al_i' + \beta_i \log N
\le \nsup{\bra'-\bra} +  \beta_i \log N.
\]
On termine en remarquant que
\dsm{
\beta_i = \frac{a}{p_k^\ro-1} \le \frac{2a}{p_k^\ro}
= \bigo{\frac{1}{(\log N)^{1+o(1)}}},
}
vu \eqref{lpk}
\hfill$\Box$

Le lemme \ref{climit2} ci-dessous étend aux voisins d'un champion
les propriétés des nombres $K$-champions énoncées dans le
lemme \ref{climit}.

\begin{lem}\label{climit2}
Avec les notations de la définition \ref{defvoisin},
soit $N$ un $K$-champion qui tend vers l'infini et $M =
2^{\al'_1\al'_1} 3^{\al'_2} \dots p_{k'}^{\al'_{k'}}$ un
voisin de $N$, et   $\bra' = (\al_1',\al_2',\dots,\al'_{k'})$,
(ou plus généralement,
$\bra' \in \Rplus^{\star k'}$ avec
$\nsup{\bra'-\bra} = \bigo{(\log N)^{\mu}}$).
Soit $\delta$, $a$,$T_0$ et $B_0$ les constantes figurant dans le
lemme \ref{climit}.  Alors on a
\begin{enumerate}
\item
\begin{equation}\label{cbeta2}
c(\bra') = a \log N   + \bigo{(\log N)^{\delta}},
\end{equation}
\item
\begin{equation}\label{tbeta2}
T(\bra') = T_0 +  \bigo{(\log N)^{\delta-1}},
\end{equation}
\item
\begin{equation}\label{bbeta2}
\BB(\bra') = B_{0} \sqrt{\log N}
\left(1+\bigo{(\log N)^{\delta-1}}\right),  
\end{equation}
\item
\begin{equation}\label{Kchampionequiv2}
K(M) =  B_{0}\sqrt\pi \sqrt{\log N}\,\, \AA(\bra')
\left(1+\bigo{(\log N)^{\delta-1}} \right).
\end{equation}
\end{enumerate}
\end{lem}

\paragraph{Démonstration de \eqref{cbeta2}}
Vu \eqref{cbeta}, il suffit de montrer que 
$\abs{c(\bra')-c(\bra)} = \bigo{(\log N)^{\delta}}$.
Le lemme \ref{lemcxy} et la définition \ref{defvoisin} donnent
\[
\abs{c(\bra') - c(\bra)} \le 2\nsup{\bra'-\bra} 
= \bigo{(\log N)^{\delta}}. 
\]

\paragraph{Démonstration de \eqref{tbeta2}}
Vu \eqref{tbeta}, il suffit
de montrer que $\abs{T(\bra')-T(\bra)} = \bigo{(\log N)^{\delta-1}}$.
Le lemme \ref{Tlipsch} donne
\[
\abs{T(\bra') - T(\bra)} \le 
3 \frac{\nsup{\bra'-\bra}}{\max(\Omega(N),\Omega(M))}\cdot
\]
Par \eqref{decNth1.2} et \eqref{bigomegaequiv} on a
$\Omega(M) \sim \Omega(N) \sim a \log N$, et avec
la définition \ref{defvoisin}, on en déduit
\[
\abs{T(\bra') - T(\bra)} = \bigo{\frac{(\log N)^\delta}{\log N}}
= \bigo{(\log N)^{\delta-1}}.
\]

\paragraph{Démonstration de \eqref{bbeta2}}
La formule \eqref{bbeta2} se déduit de la définition \eqref{defB}
de $\BB$, de \eqref{cbeta2} et de \eqref{tbeta2}.

\paragraph{Démonstration de \eqref{Kchampionequiv2}}
En choisissant $\eta=1-\delta$ dans le théorème \ref{evans}
on obtient
\begin{equation*}
K(M) = 
\sqrt{\pi}\BB(\bra')\AA(\bra')\left(1+\bigo{(\Omega(M))^{\delta-1}}\right)\\
\end{equation*}
On conclut avec \eqref{bbeta2} et \eqref{bigomegaequiv}.
\hfill$\Box$

\begin{lem}
Soit $N = 2^{\al_1} 3^{\al_2} \dots p_{k}^{\al_{k}}$ un $K$-champion
tendant vers l'infini et $\bra = (\al_1,\al_2,\dots,\al_k)$.
Soit alors $\brx \in \Rplus^{\star \ell}$ dépendant de $N$, tel que
$\nsup{\brx-\bra} = \bigo{(\log N)^{\delta}}$. 
Alors, pour tout $i$ tel que 
$p_{i}^{\ro} = o\left((\log N)^{1-\delta}\right)$ on a 
\begin{equation}\label{fklim1}
\dron{F(\brx)}{x_i} = \ro\log p_i + \bigo{p_{i}^{\ro}(\log N)^{\delta-1}}.
\end{equation}
\end{lem}

\begin{dem}
Par \eqref{defbetai}, on a $\beta_i \asymp \frac{1}{p_{i}^{\ro}}$.
Par le lemme \ref{lemcxy}, \eqref{decNth1.1}, \eqref{cbeta}
et la définition de $\beta_i$ dans \eqref{defbetai},
il vient 
\begin{eqnarray*}
\frac{c(\brx)+x_i}{x_i} &=& \frac{c(\bra) 
 + \al_i +\bigo{(\log N)^{\delta}}}{\al_i+\bigo{(\log N)^{\delta}}}\\
&=& 
\frac{a \log N + \beta_i \log N + \bigo{(\log N)^{\delta}}}{\beta_i \log N
  + \bigo{(\log N)^{\delta}}}\\
&=& \frac{(a+\beta_i)\left(1+\bigo{(\log N)^{\delta-1}} \right)}
{\beta_i\left(1+\bigo{p_i^\ro(\log N)^{\delta-1}} \right)}\\
&=&p_i^{\ro}\left(1 + \bigo{p_i^\ro(\log N)^{\delta-1}}\right),
\end{eqnarray*}
et donc, par la formule \eqref{d1F}, 
\[
\dron{F}{x_i}(\brx) = 
\log\left(p_i^{\ro}\left(1 + \bigo{p_i^\ro(\log N)^{\delta-1}}\right)\right)
= \ro\log p_i + \bigo{p_i^\ro(\log N)^{\delta-1}}.
\]
\end{dem}

Le lemme précédent précise la variation de $F(\bra)$
lorsque on multiplie un champion $N= 2^{\al_1}3^{\al_3}\dots
p_k^{\al_k}$ par un petit facteur premier. Le lemme suivant précise
la variation de $F(\bra')$ lorsque l'on multiplie
$M'=2^{\al_1'}3^{\al_2'}\dots p_{k'}^{\al_{k'}'}$, voisin d'un champion, par
un grand facteur premier.

\begin{lem}
Soit $N=2^{\al_1}3^{\al_2}\dots p_k^{\al_k}$ un $K$-champion
tendant vers l'infini et
$M' = 2^{\al'_1}3^{\al'_2}\dots p_{k'}^{\al'_{k'}}$ un voisin de $N$.
Soit $i$ un indice tel que $p_{i}^{\ro} \gg (\log N)^{1-\delta}$
et $i  \le k'+1$. Soit 
$M^{\sg} = M' p_i =  2^{\al\sg_1}3^{\al\sg_2}\dots p_{k\sg}^{\al_{k\sg}\sg}$
avec $k\sg = \max(i,k')$. En posant
$\bra\sg = (\al_1\sg,\al_2\sg,\dots,\al_{k\sg}^{\sg})$ 
et $\bra' = (\al_1',\al_2',\dots,\al_{k'}^{'})$.
Lorsque $i=k'+1$, on pose $\al'_{k'+1} = 0$\,
et $\bra' = (\al_1',\al_2',\dots,\al'_{k'},0)$.
On a alors
\begin{multline}\label{fklim2}
F(\bra\sg) - F(\bra') =  1 + \log a + \log\log N\\
-(\al_i'+1)\log(\al_i'+1) + \al_i'\log \al_i'
+  \bigo{(\log N)^{\delta-1}}
\end{multline}
en convenant que $\al_i'\log \al_i' = 0$ pour $\al_i'=0$.
\end{lem}

\noindent{\bf Démonstration}
Commençons par remarquer que $M\sg$ est aussi un voisin de $M$.
\begin{enumerate}
\item
Il résulte de $p_i^\ro \gg (\log N)^{1-\delta}$, de \eqref{decNth1.1} et de la
définition \ref{defvoisin} que l'on a
\begin{equation}\label{jzero}
\al_i' = \bigo{(\log N)^{\delta}}.
\end{equation}
Notons $c\sg = c(\bra\sg)$ et $c' = c(\bra')$.
Montrons que 
\begin{equation}\label{cmc}
c\sg-c' = \frac{1}{T_0} + \bigo{(\log N)^{\delta-1}}.
\end{equation}
On considère la fonction $G(t)$ définie par $G(t) = c(\bra' +
t(\bra\sg-\bra'))$. Elle est dérivable et, par le théorème
des accroissements finis, $c\sg - c' = G(1) - G(0) = G'(t)$, avec $t
\in ]0,1[$.  Notons $\bar\gamma = \bra' +
t(\bra\sg-\bra')$.  Par le lemme \ref{derivc}, \eqref{cbeta2},
\eqref{tbeta2} et \eqref{jzero} il vient

\begin{eqnarray*}
G'(t) &=& \dron{c(\brx)}{x_i}(\bar\gamma) = 
\frac{1}{T(\bar\gamma)} \frac{c(\bar\gamma)}{c(\bar\gamma)+ \al_i'+t}\\
&=&
\frac{1}{T_0 +\bigo{(\log N)^{\delta-1}}}
\frac{a\log N + \bigo{(\log N)^{\delta}}}{a\log N + \bigo{(\log N)^{\delta}}}
\\
&=&
\frac{1}{T_0} +\bigo{(\log N)^{\delta-1}}
\end{eqnarray*}
ce qui démontre \eqref{cmc}.
\item
Ecrivons, par la définition \ref{defF}
\begin{multline*}
F(\bra\sg)-F(\bra')
=\sum_{j=1,j \ne i}^{k\sg} \al_j'\log\frac{c\sg+\al_j'}{\al_j'}
- \sum_{j=1, j \ne i}^{k\sg} \al_j'\log\frac{c'+\al_j'}{\al_j'}\\
 + (\al_i'+1)\log\left(\frac{c\sg+\al_i'+1}{\al_i'+1}\right) 
    - \al_i'\log\left(\frac{c'+\al_i'}{\al_i'}\right)
\end{multline*}
soit
\begin{equation}\label{goodsum}
F(\bra\sg)-F(\bra') = S_1 + S_2 + \al_i'\log \al_i' - (\al_i'+1)\log(\al_i'+1)
\end{equation}
avec
\begin{equation}\label{s1def}
S_1 = 
\sum_{j=1,j\ne i}^{k\sg}\al_j'\log\left(1+\frac{c\sg-c'}{c'+\al_j'}\right)
\end{equation}
et
\begin{equation}\label{s2def}
S_2 =  (\al_i'+1)\log(c\sg+\al_i'+1) - \al_i'\log\left(c'+\al_i'\right).
\end{equation}
Par \eqref{bigomegaequiv} et \eqref{cbeta2},\,
$c'+\al_j' \asymp \log N$, ce qui entraîne par \eqref{cmc}, 
\begin{equation*}
\log\left(1+\frac{c\sg-c'}{c'+\al_j'}\right) 
= \frac{c\sg-c'}{c'+\al_j'} + \frac{\bigo{1}}{(\log N)^2}
= \frac{1}{T_0}\frac{1}{c'+\al_j'}
+ \bigo{(\log N)^{\delta-2}}\cdot
\end{equation*}
Il en résulte que 
\begin{eqnarray*}
S_1 &=& \sum_{j=1,j \ne i}^{k}\al_j' \log\frac{c\sg+\al_j'}{c'+\al_j'}\\
&=& \frac{1}{T_0} \sum_{j=1, j \ne i}^{k} \frac{\al_j'}{c'+\al_j'} +
\left(\sum_{j=1,j \ne i}^{k} \al_j' \right)\bigo{(\log
N)^{\delta-2}}\cdot
\end{eqnarray*}
Puis, en utilisant \eqref{dronc}, \eqref{tbeta2},
\eqref{bigomegaequiv} et \eqref{jzero}
\begin{equation}\label{goodsum2}
S_1 
= \frac{1}{T_0}\left(T(\bra')-\frac{\al_i'}{c'+\al_i'} \right)
+ (\Omega(M')-\al_i')\bigo{(\log N)^{\delta-2}}
= 1 + \bigo{(\log N)^{\delta-1}}\cdot
\end{equation}
\item
La définition \eqref{s2def} de $S_2$ donne
\begin{eqnarray*}
S_2 &=& \log(c\sg+\al_i'+1) + \al_i' \log\frac{c\sg+\al_i'+1}{c'+\al_i'} \\
    &=&  \log(c\sg+\al_i'+1) + \al_i' \log\left(1+\frac{c\sg-c'+1}{c'+\al_i'}\right)\cdot
\end{eqnarray*}
Avec \eqref{cbeta2}, \eqref{cmc} et \eqref{jzero} on en déduit
\begin{eqnarray*}
S_2 &=&
 \log\left(a\log N(1+O(\log N)^{\delta-1})\right) +
 \bigo{(\log N)^{\delta}}\bigo{(\log N)^{-1}} \\
&=& \log a + \log\log N + \bigo{(\log N)^{\delta-1}}\cdot
\end{eqnarray*}
Avec \eqref{goodsum} et \eqref{goodsum2} on en déduit
\eqref{fklim2}.
\hfill$\Box$
\end{enumerate}

\begin{lem}\label{Nfoisp}
  Soit $N = 2^{\al_1} 3^{\al_2} \dots p_{k}^{\al_{k}}$ un $K$-champion
  tendant vers l'infini et
  $M' = 2^{\al_1'}3^{\al_2'}\dots p_{k'}^{\al'_{k'}}$ un 
  voisin de $N$. 
\begin{enumerate}
\item
Soit $p_i$ un nombre premiet tel que $p_i =\bigo{(\log N)^{\eta}}$,
avec
\begin{equation}\label{zerodouze}
O \le \eta < \frac{1-\delta}{\ro} = 0.12192\dots
\end{equation}
Soit $u$ petit entier, $u = \bigo{(\log\log N)}$,
et $M\sg = p_i^uM'$. Alors, lorsque $N$ tend vers l'infini,
\begin{equation}\label{mulpetit}
\frac{K(M\sg)}{K(M')} =  p_i^{\ro u}
\left(1+\bigo{\frac{\log\log N}{(\log N)^{1-\delta-\eta\ro}}}\right)
\cdot
\end{equation}
\item
Si $p_i \gg (\log N)^{(1-\delta)/\ro}$ et $i \le k'+1$ alors,
lorsque $N \to \infty$, on a, en posant
$\al'_{k'+1} = 0$,
\begin{equation}\label{mulgrand}
\frac{K(M\sg)}{K(M')} = 
\frac{a\log N}{\al'_{i}+1}
\left(1+\bigo{(\log N)^{\delta-1}}\right)\cdot
\end{equation}
\end{enumerate}
\end{lem}

\begin{dem}
\begin{enumerate}
\item 
Posons $\al\sg_j = \al_j'$ pour $j\ne i$, $\al\sg_i = \al_i'+u$
et $\bar\al\sg = (\al\sg_1,\al\sg_2,\dots,\al\sg_{k'})$.
Comme $M\sg$ et $M'$ sont des voisins de $N$, par
\eqref{Kchampionequiv2}, il vient
\[
\frac{K(M\sg)}{K(M')} = 
\frac{\AA(\bar\al\sg)}{\AA(\bra')}\left(1+\bigo{(\log N)^{\delta-1}} \right)
\cdot
\] 
et, par \eqref{proxfa}, on a
\begin{equation}\label{aequiv1}
\frac{A(\bar\al\sg)}{\AA(\bra')} = \exp\left(F(\bra\sg)-F(\bra')\right)
\frac{s(\al_i')}{s(\al_i'+u)}\cdot
\end{equation}
Par \eqref{defbetai} 
et l'hypothèse $p_i = \bigo{(\log N)^{\eta}}$ on a
\[
\beta_i = \frac{a}{p_i^\ro-1} \ge \frac{a}{p_i^\ro}
\gg \frac{a}{(\log N)^{\eta\ro}}\cdot
\]
Par \eqref{decNth1.1}, \eqref{zerodouze}, et la définition
\ref{defvoisin}, pour tout 
$j$, $0 \le j \le u$, on a $\al_i'+j+1 \gg (\log N)^{1-\eta\ro}$,
puis par le lemme \ref{lemstirling},
\[
\frac{s(\al_i'+j)}{s(\al_i'+j+1)} = 
1 - \frac{1}{2(\al_i'+j+1)} +
\bigo{\left(\frac{1}{\al_i'+j+1}\right)^2}
= 1
+ \bigo{\frac{1}{(\log N)^{1-\eta\ro}}}.
\]
De $u = \bigo{\log\log N}$ il résulte alors
\[
\frac{s(\al_i')}{s(\al_i'+u)}
 = \prod_{j=0}^{u-1} \frac{s(\al_i'+j)}{s(\al_i'+j+1)}
= 1 + \bigo{\frac{\log\log N}{(\log N)^{1-\eta\ro}}}
\]
quand $N$ tend vers l'infini, et avec \eqref{aequiv1},
\begin{equation}\label{}
\frac{\AA(\bra\sg)}{\AA(\bra')} = \exp\left(F(\bra\sg)-F(\bra')\right)
\left(1+\bigo{\frac{\log\log N}{(\log N)^{1-\eta\ro}}}\right)\cdot
\end{equation}
Par le théorème des accroissements finis, il existe $t \in ]0,1[$,
tel que
\[
F(\bra\sg)-F(\bra') = u \dron{F(\brx)}{x_i}(\bra' + t(\bra\sg -\bra')).
\]
Puisque $\nsup{\bra\sg-\bra'} = \abs{u} = \bigo{\log\log N}$,
on en déduit avec \eqref{fklim1}
\begin{eqnarray*}
F(\bra\sg)-F(\bra')
&=& u \ro\log p_i  + \bigo{(\log\log N)p_i^\ro(\log N)^{\delta-1}}\\
&=& u \ro\log p_i  + \bigo{\frac{\log\log N}{(\log N)^{1-\delta-\eta\ro}}},
\end{eqnarray*}
puis \eqref{mulpetit}, à l'aide de \eqref{zerodouze}.

\item 
Lorsque $p_i$ est grand, on a, comme dans le point 1,
\[
\frac{K(M\sg)}{K(M')} = 
\frac{\AA(\bar\al\sg)}{\AA(\bra')}\left(1+\bigo{(\log N)^{\delta-1}} \right)
\cdot
\] 
et, par \eqref{fklim2}
\begin{multline*}
F(\bra\sg)-F(\bra') = 
  1 + \log a + \log\log N\\
-(\al_i'+1)\log(\al_i'+1) + \al_i'\log \al_i'
+  \bigo{(\log N)^{\delta-1}}.
\end{multline*}
De cette égalité, et de \eqref{aequiv1} (avec $u=1$), on déduit 
en utilisant \eqref{lemj1j}
\begin{eqnarray*}
\frac{\AA(\bra\sg)}{\AA(\bra')} &=& \frac{a\log N}{\al_i'+1}
e \left(\frac{\al_i'}{\al_i'+1)}\right)^{\al_i'}
\frac{s(\al_i')}{s(\al_i'+1)}\left(1+\bigo{(\log N)^{\delta-1}}\right)\\
&=& \frac{a\log N}{\al_i'+1}\left(1+\bigo{(\log N)^{\delta-1}}\right)
\cdot
\end{eqnarray*}
\end{enumerate}
\end{dem}

\begin{theorem}\label{decNth2}
Soit $N = 2^{\al_1}3^{\al_2}\dots p_k^{\al_k}$ un $K$-champion.
\begin{enumerate}
\item
Pour $N$ suffisament grand, on a $\al_k = 1$.
\item
Pour $j \ge 1$, un entier fixé. On désigne par $P_j$
le plus grand nombre premier tel que 
$P_j^j$ divise $N$ (en particulier $P_1 = p_k$).
Lorsque $N$ tend vers l'infini
\begin{equation}\label{Pj}
P_j \sim \left(\frac{a\log N}{j}\right)^{\frac{1}{\ro}}
\sim j^{-1/\ro} P_1.
\end{equation}
\item
Lorsque $N$ tend vers l'infini,
\begin{equation}\label{omequiv}
\omega(N) \sim \ro a^{1/\ro} \frac{(\log N)^{1/\ro}}{\log\log N}\cdot
\end{equation}
\end{enumerate}
\end{theorem}

\begin{rem}
Le calcul des $761$ premiers nombres $K$-champions
laisse penser que le nombres de champions tels que
$\al_k > 1$ est $111$, et que le  plus grand d'entre
eux est le $390\iem$ champion
\[
N_{390} = 2^{28} 3^{10} 5^{4} 7^{2} = 485\;432\;135\;516\;160\;000.
\]
\end{rem}

\paragraph{Démonstartion du point 1}
On suppose $\al_k \ge 2$ et on
pose $M=N\frac{p_{k+1}p_{k+2}}{2 p_{k-1}p_k}$.  Lorsque $N$ tend
vers l'infini, $k$ tend vers l'infini.  Pour $N$ suffisament grand on
a donc $M < N$, et puisque $N$ est un champion $K(M) < K(N)$.  On pose
$M^{(0)} = N$, $M^{(1)} = \frac{N}{p_k}$, $M^{(2)} =
\frac{N}{p_{k-1}^{}p_k}$, $M^{(3)} = \frac{Np_{k+1}}{p_{k-1}^{}p_k}$,
$M^{(4)} = \frac{Np_{k+1}p_{k+2}}{p_{k-1}^{}p_k}$ et $M^{(5)} = M$.
Ces nombres sont des voisins de $N$, et, par \eqref{thomega2}, on
a $p_{k-1} \gg (\log N)^{(1-\delta)/\ro}$.
En appliquant 5 fois le lemme \ref{Nfoisp}, lorsque $N$ tend vers
l'infini, on a
\begin{equation}\label{cinq}
1 > \frac{K(M)}{K(N)}
= \prod_{i=0}^{4} \frac{K(M^{(i+1)})}{K(M^{(i)})} 
 \sim \frac{\al_{k-1}\al_{k}}{2^\ro\times 1\times 1}\cdot
\end{equation}
On a donc pour $N$ assez grand 
$
(\al_{k-1}\al_{k})/2^\ro \lesssim 1,
$
et, a fortiori,
\[
\al_k^2 \le \al_{k-1}\al_k \lesssim 2^\ro = 3.314.
\]
Puisque $\al_k$ est un entier, il est plus petit que $2$, il
y a contradiction.

\paragraph{Démonstration de \eqref{Pj}}
Soit $i$, dépendant de $N$, le plus grand indice tel que $\al_i \ge
j$. Cet entier $i$ est caractérisé par $\al_i \ge j >
\al_{i+1}$. Notons $j' = \al_i$ et $j\sg = \al_{i+1}$.  D'après
\eqref{decNth1.1}, il existe une constante $C$ telle que 
$j > \al_{i+1} \ge \beta_{i+1}\log N - C (\log N)^{\delta}
\ge \frac{a}{p_{i+1}^{\delta}} \log N - C (\log N)^{\delta}$.
Cela implique $p_{i+1}\gg (\log N)^{(1-\delta)/\ro}$, et donc
aussi $p_{i}\gg (\log N)^{(1-\delta)/\ro}$.

Le nombre $\log 3/\log 2$ étant irrationnel, si
l'on ordonne l'ensemble des $2^u 3^v$ en
une suite croissante $(x_n)_{n \ge 1}$, on a $\lim x_{n+1}/x_n =
1$. En conséquence, pour tout $\veps > 0$, il existe $N_{\veps}$ tel
que, pour tout champion $N \ge N_{\veps}$, il existe des entiers $u',
v',u\sg,v\sg \ge 0$ avec
\begin{equation}\label{uv1}
(1-\veps)p_i \le 2^{u'} 3^{v'} < p_i < p_{i+1} < 2^{u\sg}3^{v\sg}
< p_{i+1}(1+\veps) < p_i(1+2\veps).
\end{equation}
De l'encadrement ci-dessus, de $p_i \le p_k$ et de \eqref{lpk} il
résulte que
\begin{equation}\label{}
u', v', u\sg, v\sg = \bigo{\log\log N}.
\end{equation}
Soit
\begin{equation}\label{defM1}
M' = N 2^{u'} 3^{v'}/p_i < N
\quad\text{ et }\quad
M\sg = N\frac{p_{i+1}}{2^{u\sg} 3^{v\sg}} < N.
\end{equation}
Comme $N$ est un $K$-champion on a $K(M') < K(N)$
et $K(M\sg) < K(N)$.
Par le lemme \ref{Nfoisp} on a, comme en \eqref{cinq}
\[
\frac{K(M')}{K(N)} \sim \frac{2^{\ro u'} 3^{\ro v'}j'}{a \log N}
\quad\text{ et }\quad
\frac{K(M\sg)}{K(N)} 
\sim \frac{a\log N}{2^{\ro u\sg} 3^{\ro v\sg}(j\sg+1)}\cdot
\]
Et donc, puisque $2^{u'} 3^{v'} \ge (1-\veps) p_i$,
\[
(1-\veps)^{\ro} p_i^{\ro} \frac{j'}{a\log N}\lesssim 
\frac{K(M')}{K(N)} < 1.
\]
soit
\[
P_j = p_i \lesssim \frac{1}{1-\veps}\left(\frac{a\log N}{j'}
\right)^{\frac{1}{\ro}}
\le
\frac{1}{1-\veps}
\left(\frac{a\log N}{j}\right)^{\frac{1}{\ro}}\cdot
\]
De même, la majoration $K(M\sg)/K(N) < 1$ avec
$2^{u\sg}3^{v\sg} \le p_i(1+2\veps)$ donne
\[
P_j = p_i \gtrsim \frac{1}{(1+2\veps)}
\left(\frac{a\log N}{j\sg+1}\right)^{\frac{1}{\ro}} \ge
\frac{1}{(1+2\veps)}\left(\frac{a\log N}{j}\right)^{\frac{1}{\ro}}\cdot
\]

\paragraph{Démonstration de \eqref{omequiv}}
 
Choisissant $j=1$ dans \eqref{Pj} on obtient
\[
p_k = P_1 \sim \left(a\log N\right)^{1/\ro},
\]
puis, par le théorème des nombres premiers,
$
k\log k \sim \left(a\log N\right)^{1/\ro}
$
et
$
\log k \sim \frac{1}{\ro}\log\log N
$
et cela donne\eqref{omequiv}.
\hfill$\Box$.

\begin{rem}
En utilisant les résultats connus sur les approximations 
diophantiennes de $\log 3/\log 2$ (cf. \cite{rhin}), on
peut rendre plus précises les inégalités \eqref{uv1} et 
obtenir un terme de reste pour les formules \eqref{Pj}
et \eqref{omequiv}.
\end{rem}

\subsection{Estimation de $Q(X)$}

Soit $Q(X)$ le nombre de $K$-champions au plus égaux à $X$. Par
les propriétés \eqref{factorchamp} et \eqref{facteur2},
nous avons, comme en
\cite{Her-Nic}, paragraphe 6.4
\[
\log X \ll  Q(X) \le \exp
\left((1+o(1))\frac{2\pi}{\sqrt 3}\sqrt{\frac{\log X}{\log\log X}}\right)
\]
et l'on peut montrer
\begin{theorem}
\begin{enumerate}
\item
Pour $X$ assez grand, on a
\[
Q(X) \ge (\log X)^{1.07}.
\]
\item
Lorsque $X$ tend vers l'infini
\[
\log Q(X) = \bigo{(\log X)^{\delta/2}}
\]
où $\delta = 0.788\dots$ a été défini en \eqref{defdelta}.
\end{enumerate}
\end{theorem}

\paragraph{Démonstration de 1}
La preuve est très voisine de celle de la proposition
7.1 de \cite{Her-Nic}. Soit $N$ un nombre $K$-champion
assez grand et $\eta = 3/40= 0.075$. Notons que cela
entraine $1-\delta-\eta\ro = 0.0811\ldots > \eta + 0.006.$
On montre d'abord que
le nombre $K$-champion $N'$ suivant $N$ vérifie
\begin{equation}\label{Q1}
N' \le N\left(1+ \frac{2\eta\log\log N}{(\log N)^{\eta}} \right)
\end{equation}
ce qui, comme dans \cite{Her-Nic}, entraîne le point 1.
Pour prouver \eqref{Q1} on construit deux nombres premiers consécutifs
$p_r$ et $p_{r+1}$ vérifiant (cf. \cite{Her-Nic})
\begin{equation}\label{Q2}
(\log N)^{\eta} \le p_r < p_{r}+2 \le p_{r+1} < 2(\log N)^{\eta}
\end{equation}
et\begin{equation}\label{Q3}
p_r < p_{r}+2 \le p_{r+1} \le p_r + 2\eta\log \log N.
\end{equation}
On pose $M = \frac{p_{r+1}}{p_r} N$;
on s'assure par \eqref{decNth1.1} que l'exposant $\al_r$ de $p_r$ dans
$N$ tend vers l'infini, et que $N/p_r$ et $M$ sont
des voisins de $N$, puis, par l'égalité
\eqref{mulpetit} du lemme
\ref{Nfoisp}
\begin{equation}\label{Q4}
\frac{K(M)}{K(N)} = \left(\frac{p_{r+1}}{p_r}\right)^{\ro}
\left(1+ \bigo{\frac{\log\log N}{(\log N)^{1-\delta-\eta\ro}}}\right)
=
\left(\frac{p_{r+1}}{p_r}\right)^\ro
\left(1+\frac{\bigo{1}}{(\log N)^{\eta+0.006}}\right)
\end{equation}
Mais par \eqref{Q2},
\[
\left(\frac{p_{r+1}}{p_r}\right)^{\ro} \ge \left(1+\frac{2}{p_r}\right)^{\ro}
\ge \left(1+\frac{1}{(\log N)^{\eta}}\right)^{\ro}
\ge 1 + \frac{\ro}{(\log N)^{\eta}}
\]
ce qui, avec \eqref{Q4}, montre que $K(M) > K(N)$. Cela
implique $N' \le M$ et, par \eqref{Q3}, 
\[
N' \le M = N \frac{p_{r+1}}{p_r} \le
N\left(1+\frac{2\eta\log\log N}{p_r} \right)
\]
qui, avec \eqref{Q2}, démontre \eqref{Q1}.

\paragraph{Démonstration de 2}

Soit $N = 2^{\al_1}3^{\al_2}\dots p_k^{\al_k}$, et $J$ un entier.
On pose
\[
T_J = T_J(N) = \sum_{i=J+1}^{k} \al_i.
\]
Il résulte de \eqref{majfort} que 
$\sum_{i=J+1}^{k} \abs{\al_i-\beta_i\log N}= \bigo{(\log
N)^{\delta}}$, et cela implique 
\[
T_J = \sum_{i=J+1}^{k} \beta_i \log N +  \bigo{(\log N)^{\delta}}
= \left(\sum_{i=J+1}^{k} \frac{a}{p_i^{\ro}-1}\right)\log N
+ \bigo{(\log N)^{\delta}}.
\]
Mais on a
\[
\sum_{i=J+1}^{k} \frac{a}{p_i^{\ro}-1}
\le \sum_{i=J+1}^{k} \frac{2a}{p_i^{\ro}}
\le \sum_{i=J+1}^{k} \frac{2a}{i^{\ro}}
\le \int_{J}^{\infty} \frac{2a\,dt}{t^{\ro}}
= \frac{2a}{(\ro-1)J^{\ro-1}}
\]
ce qui entraîne
\begin{equation}\label{Q5}
T_J = \bigo{\frac{\log N}{J^{\ro-1}}}
+ \bigo{(\log N)^{\delta}}.
\end{equation}
On choisit
\[
J = \intpart{(\log X)^{\gamma}}
\]
avec 
\[
\gamma = \frac{1-\delta}{\ro-1} = \frac{1}{2\ro}
= 0.289\dots
\]
Alors, pour tous les nombres $K$-champions $N \le X$,
on a par \eqref{Q5}
\[
T_J(N) = \bigo{(\log X)^{\delta}}.
\]
La preuve du point 2 se termine alors comme dans 
\cite{Her-Nic}.
$\Box$.

\begin{figure}\label{tablechampions}
\[
\begin{array}{|r||r|r|l|l|}
\hline
i & N_i & K(N_i) &\ \ \ \bra_i &\ \  K(N_i) = \\[1mm]
\hline
1&1&1&[]&1\\
2&4&2&[2]&2\\
3&6&3&[1,1]&3\\
4&8&4&[3]&2^2\\
5&12&8&[2,1]&2^3\\
6&24&20&[3,1]&2^3 \times 5\\
7&36&26&[2,2]&2 \times 13\\
8&48&48&[4,1]&2^3 \times 3\\
9&72&76&[3,2]&2^2 \times 19\\
10&96&112&[5,1]&2^4 \times 7\\
11&120&132&[3,1,1]&2^2 \times 3 \times 11\\
12&144&208&[4,2]&2^4 \times 13\\
13&192&256&[6,1]&2^8\\
14&240&368&[4,1,1]&2^4 \times 23\\
15&288&544&[5,2]&2^5 \times 17\\
16&360&604&[3,2,1]&2^2\times 151\\
17&432&768&[4,3]&2^8 \times 3\\
18&480&976&[5,1,1]&2^4 \times 61\\
19&576&1376&[6,2]&2^5 \times 43\\
20&720&1888&[4,2,1]&2^5 \times 59\\
21&864&2208&[5,3]&2^5 \times 3 \times 23\\
22&960&2496&[6,1,1]&2^6 \times 3 \times 13\\
23&1152&3392&[7,2]&2^6 \times 53\\
24&1440&5536&[5,2,1]&2^5 \times 173\\
25&1728&6080&[6,3]&2^6 \times 5 \times 19\\
26&1920&6208&[7,1,1]&2^6 \times 97\\
27&2160&7968&[4,3,1]&2^5 \times 3 \times 83\\
28&2304&8192&[8,2]&2^{13}\\
29&2880&15488&[6,2,1]& 2^7 \times 11^2\\
30&3456&16192&[7,3]&2^6 \times 11 \times 23\\
31&4320&25440&[5,3,1]&2^5 \times 3 \times 53\\
32&5760&41792&[7,2,1]&2^6 \times 653\\
33&6912&41984&[8,3]&2^{10}\times 41\\
34&8640&76864&[6,3,1]&2^6 \times 1201\\
35&11520&109568&[8,2,1]&2^{10}\times 107\\
36&17280&222528&[7,3,1]&2^6 \times 3 \times 19 \times 61\\
37&23040&280576&[9,2,1]&2^{11}\times 137\\
38&25920&331776&[6,4,1]&2^{12}\times 3^4\\
39&30240&333984&[5,3,1,1]&2^5 \times 3 \times 7^2 \times 71\\
40&34560&622592&[8,3,1]&2^{15}\times 19\\
\hline
\end{array}
\]
\caption{Table des 40 premiers $K$-champions}
\end{figure}

\subsection{Table des champions}

Nous avons vu \eqref{factorchamp} que les nombres $K$-champions sont
de la forme
\[
N = 2^{\al_1}3^{\al_2}\dots p_k^{\al_k},
\quad\text{ avec }\quad
\al_1 \ge \al_2 \ge \cdots \ge \al_k \ge 1.
\] Pour énumérer tous les $K$ champions inférieurs à $X$ on construit,
par \og backtracking \fg\ tous les nombres $N$ de cette forme
jusqu'à $X$.  Le bactracking n'énumère pas ces nombres par
ordre croissant.  Chaque fois qu'un tel $N$ est obtenu, la valeur
$K(N)$ est calculée en utilisant la formule \eqref{KP} et le couple
$(N,K(N))$ est mémorisé.  Une fois cette construction
terminée, on range ces couples par valeurs croissantes de $N$, et
on élimine ceux qui ne sont pas champions.  Nous avons calculé
tous les champions jusqu'à $X =
557\,940\,830\,126\,698\,960\,967\,415\,390$.  Nous avons obtenu
$340\,884$ nombres de la forme \eqref{factorchamp}.  Parmi ceux-ci
$761$ sont des champions.  La figure 2 donne les 40 premiers
champions.
\footnote{ 
On trouvera une table plus complète sur l'une ou l'autre
des deux pages personnelles
\centerline{\tt http\!\!\!://math.univ-lyon1.fr/$\,\,\widetilde{}\,\,$deleglise}
\centerline{\tt http\!\!\!://math.univ-lyon1.fr/$\,\,\widetilde{}\,\,$nicolas\ \ }
}

\subsection{Problèmes ouverts}
Les problèmes listés à la fin de l'article \cite{Her-Nic}
au sujet de la fonction $K_{\cal P}$ (définie en \eqref{KP})
peuvent se poser à peu près dans les mêmes termes
pour la fonction de Kalm\'ar $K$.

\newcommand{\bibit}[5]{\bibitem{#1} {\sc #2} #3 {\sl #4} #5}

\vfill

\hrule
\medskip
\noindent
\begin{small}
\begin{minipage}[t]{5.3cm}
deleglise{\at}math.univ-lyon1.fr\\[0.5mm]
jlnicola{\at}in2p3.fr
\end{minipage}
\begin{minipage}[t]{7cm}
Université de Lyon,
Université Lyon1, CNRS.\\
UMR 5208, Institut Camille Jordan\\
Bât. Jean Braconnier\\
21 Avenue Claude Bernard\\
F-69622 Villeurbanne cedex, \ France.
\end{minipage}
\medskip

\noindent
\begin{minipage}[t]{5.3cm}
mhernane@usthb.dz
\end{minipage}
\begin{minipage}[t]{7cm}
Université des Sciences et de la Technologie\\
Institut de Mathématiques.\\
BP 32, USTHB, 16123 Bab-Ezzouar\\
Alger, \ Algérie.
\end{minipage}
\end{small}

\end{document}